\documentclass[draftclsnofoot,onecolumn]{IEEEtran}

\usepackage{ifpdf,graphicx,amsmath,amsbsy,amssymb,cite,citesort,epsf,url,color,epstopdf}
\usepackage{amscd,amsfonts,amsgen,amsopn,amstext,amsthm,amsxtra}
\usepackage{graphicx}
\usepackage{verbatim}
\usepackage{algorithm}
\usepackage{algorithmic}
\theoremstyle{plain}
\newtheorem{thm}{Theorem}
\newtheorem{lem}{Lemma}
\newtheorem{cor}{Corollary}

\theoremstyle{definition}
\newtheorem{defi}{Definition}
\newtheorem{rem}{Remark}

\def \real { \mathbb{R} }
\def \y {y} 

\def \thetahat
\def \thresh {{\mathcal{T}}}
\def \thetahat {\widehat{\theta}}
\def \betahat {\widehat{\beta}}

\def \pinv {^\dag}

\def \prob {\mathbb{P}}

\def \Hhat {\widetilde{H}}
\def \ident {\mathrm{Id}}

\newcommand{\cl}[1]{{\mathcal{#1}}}
\newcommand{\Eq}{{\mathbb{E}_q}}
\newcommand{\E}{{\mathbb{E}}}

\newcommand{\R}{{\mathbb{R}}}

\newcommand{\sign}{{\mathrm{sign}}}
\newcommand{\signbar}{{\overline{\mathrm{sign}}}}

\begin{document}
\title{\huge Conditioning of Random Block Subdictionaries with Applications to Block-Sparse Recovery and Regression}
\author{Waheed U. Bajwa,~\IEEEmembership{Senior~Member,~IEEE,}
Marco F. Duarte,~\IEEEmembership{Senior~Member,~IEEE,}\\
and~Robert Calderbank,~\IEEEmembership{Fellow,~IEEE}
\thanks{The first two authors contributed equally to the paper and are listed in alphabetical order. Portions of this work have previously appeared in technical reports~\cite{BlockSubdictionariesTR,GroupTR}, and at the Int. Conf. Sampling Theory and Applications (SAMPTA), 2011~\cite{DBC11} and Int. Conf. Artificial Intelligence and Statistics (AISTATS), 2014~\cite{Bajwa.etal.Conf2014}. This work was supported by the NSF under grant CCF-1218942, by the ARO under grant W911NF-14-1-0295, by the ONR under grant N00014-08-1-1110 and by the AFOSR under grants FA9550-09-1-0422, FA9550-09-1-0643, and FA9550-05-0443. MFD was also supported by NSF Supplemental Funding DMS-0439872 to UCLA-IPAM, P.I.\ R.\ Caflisch.}%
\thanks{W. U. Bajwa is with the Department of Electrical and Computer Engineering, Rutgers University, Piscataway, NJ 08854. E-mail: {\tt waheed.bajwa@rutgers.edu}}%
\thanks{M. F. Duarte is with the Department of Electrical and Computer Engineering, University of Massachusetts, Amherst, MA 01003. E-mail: {\tt mduarte@ecs.umass.edu}}%
\thanks{R. Calderbank is with the Departments of Computer Science, Electrical and Computer Engineering, and Mathematics at Duke University, Durham, NC 27708. E-mail: {\tt robert.calderbank@duke.edu}}%
}

\maketitle

\begin{abstract}
The linear model, in which a set of observations is assumed to be given by a linear combination of columns of a matrix (often termed a dictionary), has long been the mainstay of the statistics and signal processing literature. One particular challenge for inference under linear models is understanding the conditions on the dictionary under which reliable inference is possible. This challenge has attracted renewed attention in recent years since many modern inference problems (e.g., high-dimensional statistics, compressed sensing) deal with the ``underdetermined'' setting, in which the number of observations is much smaller than the number of columns in the dictionary. This paper makes several contributions for this setting when the set of observations is given by a linear combination of a small number of \emph{groups of columns} of the dictionary, termed the ``block-sparse'' case. First, it specifies conditions on the dictionary under which most block submatrices of the dictionary (often termed block subdictionaries) are well conditioned. This result is fundamentally different from prior work on block-sparse inference because ($i$) it provides conditions that can be explicitly computed in polynomial time, ($ii$) the given conditions translate into near-optimal scaling of the number of columns of the block subdictionaries as a function of the number of observations for a large class of dictionaries, and ($iii$) it suggests that the spectral norm and a related measure of quadratic-mean block coherence of the dictionary (rather than the worst-case column/block coherences) fundamentally limit the scaling of dimensions of the well-conditioned block subdictionaries. Second, in order to help understand the significance of this result in the context of block-sparse inference, this paper investigates the problems of block-sparse recovery and block-sparse regression in underdetermined settings. In both of these problems, this paper utilizes its result concerning conditioning of block subdictionaries and establishes that near-optimal block-sparse recovery and block-sparse regression is possible for certain dictionaries as long as the dictionary satisfies easily computable conditions and the coefficients describing the linear combination of groups of columns can be modeled through a mild statistical prior. Third, the paper reports numerical experiments that highlight the effects of different dictionary measures in block-sparse inference problems.
\end{abstract}

\begin{IEEEkeywords}
Block-sparse inference, block-sparse recovery, block-sparse regression, compressed sensing, group lasso, group sparsity, high-dimensional statistics, multiple measurement vectors, random block subdictionaries
\end{IEEEkeywords}

\section{Introduction}
Consider the classical linear forward model $y = \Phi\beta$, which relates a parameter vector $\beta \in \R^p$ to an observation vector $y \in \R^n$ through a linear transformation (henceforth referred to as a \emph{dictionary}) $\Phi \in \R^{n \times p}$. This forward model, despite its apparent simplicity, provides a reasonable mathematical approximation of reality in a surprisingly large number of application areas and scientific disciplines \cite{Kay.Book1993,Kay.Book1998,Rencher.Schaalje.Book2008}. While the operational significance of this linear (forward) model varies from one application to another, the fundamental purpose of it in all applications stays the same: \emph{given knowledge of $y$ and $\Phi$, make an inference about $\beta$.} However, before one attempts to solve an inference problem using the linear model, it is important to understand the conditions under which doing so is even feasible. For instance, inferring anything about $\beta$ will be a moot point if the nullspace of $\Phi$ were to contain $\beta$. Thus, a large part of the literature on linear models is devoted to characterizing conditions on $\Phi$ and $\beta$ that facilitate reliable inference.

Classical literature on inference using linear models proceeds under the assumption that the number of observations $n$ equals or exceeds the number of parameters $p$. In this setting, conditions such as $\Phi$ being full column rank or $\Phi \Phi^*$ being well conditioned---both of which can be explicitly verified---are common in the inference literature \cite{Kay.Book1993,Knight.Fu.AS2000,Zou.JASA2006}. In contrast, there has recently been a growing interest to study inference under linear models when $n$ is much smaller than $p$. This setting is the hallmark of high-dimensional statistics \cite{Donoho.Conf2000}, arises frequently in many application areas \cite{Daubechies-Defrise-DeMol-04}, and forms the cornerstone of the philosophy behind compressed sensing \cite{DonohoCS,CandesCS}. It of course follows from simple linear algebra that inferring about every possible $\beta$ from $\y = \Phi\beta$ is impossible in this setting; instead, the high-dimensional inference literature commonly operates under the assumption that $\beta$ has only a few nonzero parameters---typically on the order of $n$---and characterizes corresponding conditions on $\Phi$ for reliable inference. Some notable conditions in this regard include the spark \cite{Donoho.Elad.PNAS2003}, the restricted isometry property \cite{CandesRIP}, the irrepresentable condition \cite{Zhao.Yu.JMLR2006} (and its variant, the incoherence condition \cite{Wainwright.ITIT2009a}), the restricted eigenvalue assumption \cite{Bickel.etal.AS2009}, and the nullspace property \cite{cohen}. While these and other conditions in the literature differ from each other in one way or the other, they all share one simple fact: \emph{requiring that $\Phi$ satisfies one of these conditions implies that one or more column submatrices (subdictionaries) of $\Phi$ must be full column rank and/or well conditioned.} Unfortunately, explicitly verifying that $\Phi$ satisfies one of these properties is computationally daunting (NP-hard in some cases \cite{Bandeira.etal.unpublished2012}), while indirect means of verifying these conditions provide rather pessimistic bounds on the dimensions of subdictionaries of $\Phi$ that are well conditioned \cite{Bajwa.Pezeshki.FF2012}.

In a recent series of influential papers, several researchers have managed to circumvent the pessimistic bounds associated with verifiable conditions on $\Phi$ for high-dimensional inference by resorting to an \emph{average-case analysis} \cite{Candes04C,CandesQRUP,Tropp.ACHA2008,TroppCRAS,Kuppinger.etal.ITIT2012,Gurevich.Hadani.unpublished2009}. Representative work by Tropp \cite{Tropp.ACHA2008,TroppCRAS}, for instance, shows that \emph{most} subdictionaries of $\Phi$ having unit $\ell_2$-norm columns are guaranteed to be well conditioned when the number of columns in the subdictionary is proportional to $p/(\|\Phi\|^2_2 \log{p})$---provided that correlations between the columns of $\Phi$ do not exceed a certain threshold, a condition readily verifiable in polynomial time. In particular, these results imply that if $\Phi$ is a unit norm tight frame \cite{Bajwa.Pezeshki.FF2012}, corresponding to $\|\Phi\|_2^2 = p/n$, then it can be \emph{explicitly verified} that most subdictionaries of $\Phi$ of dimension $n \times O(n/\log{p})$ are well conditioned.\footnote{Recall Landau's notation: $f(x) = O\big(g(x)\big)$ if there exist some $c_0, x_0$ such that $f(x) \leq c_0 g(x)$ for all $x \geq x_0$.} The biggest advantage of such {\em average-case analysis} results for the conditioning of subdictionaries of $\Phi$ lies in their ability to facilitate tighter verifiable conditions for inference under the linear model using an arbitrary (random or deterministic) dictionary $\Phi$. Several works in this regard have been able to leverage the results of \cite{Tropp.ACHA2008,TroppCRAS} to provide tighter verifiable conditions for average-case sparse recovery \cite{Herman.Strohmer.ITSP2009,Pfander.Rauhut.JFAA2010} (i.e., obtaining $\beta$ from $y=\Phi\beta$ \emph{\`{a} la} compressed sensing \cite{DonohoCS,CandesCS}), average-case model selection \cite{CandesPlan} (i.e., estimating locations of the nonzero entries of $\beta$ from $y=\Phi\beta + \text{noise}$), and average-case linear regression \cite{CandesPlan} (i.e., estimating $\Phi \beta$ from $y=\Phi\beta + \text{noise}$).

\subsection{Our Contributions}\label{ssec:contributions}
Our focus in this paper is on inference under the linear model in the ``$n$ smaller than $p$'' setting, in the case when $\beta$ not only has a few nonzero parameters, but also its nonzero parameters exhibit a certain \emph{block} (or \emph{group}) structure. Specifically, we have $\beta = [\beta_1^*~\beta_2^*~\ldots~\beta_r^*]^*$ with $\beta_i \in \R^m$ for $m, r \in \mathbb{Z}_+$, $p = rm$, and only $k \ll r$ of the $\beta_i$'s are nonzero (sub)vectors. Such setups are often referred to as \emph{block sparse} (or \emph{group sparse}) and arise in various contexts in a number of inference problems \cite{YuanLin,BachGroup,MishaliEldar,Bolstad.etal.ITS2011,Duarte.Eldar.ITSP2011}. The most fundamental challenge for inference in this block-sparse setting then becomes specifying conditions under which one or more \emph{block subdictionaries} of $\Phi$ are full column rank and/or well conditioned. A number of researchers have made substantial progress in this regard recently, reporting conditions on $\Phi$ in the block setting that mirror many of the ones reported in \cite{Donoho.Elad.PNAS2003,CandesRIP,Zhao.Yu.JMLR2006,Wainwright.ITIT2009a,Bickel.etal.AS2009,cohen} for the classical setup; see, e.g., \cite{YuanLin,BachGroup,Cotter_Sparse,TGS06,T06,GRSV08,nardi:ejs08,MeierGroup,LiuZhang,Stojnic.etal.ITSP2009,EM09,ER10,BCDH10,Stojnic2,EKB10,ZhangGroup,fusion,FangLi,BenHaim,Kim,DE10,ObozinskiUnion,Kolar,Fang,Lounici.etal.AS2011,Lee.etal.ITIT2012,Elhamifar.Vidal.ITSP2012,DWBSB13}. However, just like in the classical setup, verifying that $\Phi$ satisfies one of these properties in the block setting ends up being either computationally intractable or results in rather pessimistic bounds on the dimensions of block subdictionaries of $\Phi$ that are well conditioned. In contrast to these works, and in much the same way \cite{Candes04C,CandesQRUP,Tropp.ACHA2008,TroppCRAS,Kuppinger.etal.ITIT2012,Gurevich.Hadani.unpublished2009} reasoned in the classical case, we are interested in overcoming the pessimistic bounds associated with verifiable conditions on $\Phi$ for high-dimensional inference in the block-sparse setting by resorting to an average-case analysis.

Our first main contribution in this regard is a generalization of \cite{Tropp.ACHA2008,TroppCRAS} that establishes that \emph{most block subdictionaries} of $\Phi$ having unit $\ell_2$-norm columns are guaranteed to be well conditioned with the number of \emph{blocks} in the subdictionary proportional to $\min \{1/(\overline{\mu}_B^2 \log{p}), r/(\|\Phi\|^2_2 \log{p})\}$ provided that $\Phi$ satisfies a polynomial-time verifiable condition that we term the \emph{block incoherence condition}. Here, $\overline{\mu}_B$ denotes the quadratic-mean block coherence of the dictionary $\Phi$, which will be formally defined later. In particular, these results also imply that if $\Phi$ is a unit norm tight frame with $\overline{\mu}_B = O(\sqrt{\frac{m}{n}})$ then it can be explicitly verified that most \emph{block subdictionaries} of $\Phi$ of dimension $n \times O(n/\log{p})$ are well conditioned.

While our ability to guarantee that most block subdictionaries of a dictionary that satisfies the block incoherence condition are well conditioned makes us optimistic about the use of such dictionaries in inference problems, there remains an analytical gap in going from conditioning of block subdictionaries to performance of inference tasks. Our second main contribution in this regard is the application of results concerning the conditioning of block subdictionaries to provide tighter verifiable conditions for average-case block-sparse recovery (i.e., obtaining $\beta$ from $y=\Phi\beta$ with $\beta$ being block sparse) and average-case block-sparse regression (i.e., estimating $\Phi \beta$ from $y=\Phi\beta + \text{noise}$ with $\beta$ being block sparse).

Last, but not least, we carry out a series of numerical experiments to highlight an aspect of inference under the linear model that is rarely discussed in the related literature: \emph{the spectral norm (and the related quadratic-mean block coherence measure) of the dictionary $\Phi$ influence the inference performance much more than any of its other measures}. Specifically, our results show that performances of block-sparse recovery and regression are inversely proportional to $\|\Phi\|_2^2$ and tend to be independent of correlations between individual columns of $\Phi$ for the most part---an outcome that also hints at the possible (orderwise) tightness of our results concerning the conditioning of block subdictionaries.

\subsection{Notational Convention and Organization}
The following notation will be used throughout the rest of this paper. We use uppercase and lowercase Roman/Greek letters for matrices and vectors/scalars, respectively. Given a vector $v$, we use $\|v\|_q$ and $v^*$ to denote the usual $\ell_q$ norm and conjugate transpose of $v$, respectively. We define the \emph{scalar} sign operator for $x \in \mathbb{R}$ as $\sign(x) := x/|x|$, while we use $\sign(v)$ for a vector $v$ to denote entry-wise sign operation. In addition, we define the \emph{vector} sign operator for a vector $v$ as $\signbar(v) := v/\|v\|_2$, which returns the unit-norm vector pointing in the direction of $v$. Given two vectors $u$ and $v$, we define the inner product between them as $\langle u,v\rangle := \sum_{i} u_i v_i$. Given a matrix $A$, we use $\|A\|_2$ and $A^*$ to denote the spectral norm ($\sigma_{\max}(A)$) and the adjoint operator of $A$, respectively. In addition, assuming $A$ has unit $\ell_2$-norm columns and using $A_i$ to denote the $i^{th}$ column of $A$, the \emph{coherence} of $A$ is defined as $\mu(A) := \max_{i,j:i\not=j} |\langle A_i,A_j \rangle|$. Given a set $\cl{S}$, we use $A_\cl{S}$ (resp. $v_\cl{S}$) to denote the submatrix (resp. subvector) obtained by retaining the columns of $A$ (resp. entries of $v$) corresponding to the indices in $\cl{S}$. Given a random variable $R$, we use $\Eq[R]$ to denote $\big(\mathbb{E}[R^q]\big)^{1/q}$. Given a function $f(R)$, sometimes we also use $\Eq_{,R}[f(R)]$ to explicitly specify the random variable over which the expectation is being taken. Finally, $\ident$ denotes the identity operator and $\otimes$ denotes a Kronecker product.

The rest of this paper is organized as follows. Section~\ref{sec:result} presents the main result of this paper concerning the conditioning of block subdictionaries. Section~\ref{sec:grouprecovery} leverages the result of Section~\ref{sec:result} and presents an average-case analysis of convex optimization-based block-sparse recovery from noiseless measurements, along with some discussion and numerical experiments. Section~\ref{sec:grouplasso} makes use of the result of Section~\ref{sec:result} to present an average-case analysis of block-sparse regression and the associated numerical experiments. Finally, some concluding remarks are provided in Section~\ref{sec:conc}. For the sake of clarity of exposition, we relegate the proofs of most of the lemmas and theorems to several appendices.

\section{Conditioning of Random Block Subdictionaries}
\label{sec:result} In this section, we state and discuss the main result of this paper concerning the conditioning of block subdictionaries of the $n \times p$ dictionary $\Phi$. Here, and in the following, it is assumed that $\Phi$ has a block structure that comprises $r = p/m$ blocks of dimensions $n \times m$ each; in particular, we can write without loss of generality that $\Phi = [\Phi_1~\Phi_2~\ldots~\Phi_r]$, where each block $\Phi_i = [\phi_{i,1}~\ldots~\phi_{i,m}]$ is an $n \times m$ matrix. We also assume throughout this paper that the columns of $\Phi$ are normalized: $\|\phi_{i,j}\|_2 = 1$ for all $i=1,\dots,r$, $j=1,\dots,m$. The problem we are interested in addressing in this section is the following. Let $\cl{S} \subset \{1,\dots,r\}$ with $|\cl{S}| = k$ and define an $n \times km$ block subdictionary $X = [\Phi_i : i \in \cl{S}]$. Then what are the conditions on $\Phi$ that will guarantee that the singular values of $X$ concentrate around unity? Since addressing this question for an \emph{arbitrary} subset $\cl{S}$ is known to lead to either nonverifiable conditions or pessimistic bounds on $k$ (cf.~Section~\ref{ssec:contributions}), our focus here is on a subset $\cl{S}$ that is drawn uniformly at random from all $\binom{r}{k}$ possible $k$-subsets of $\{1,\dots,r\}$.

\subsection{Main Result}
Our main result concerning the conditioning of random block subdictionaries involves the use of three different measures of \emph{block coherence} of the dictionary.
\begin{defi}[Block Coherences]
The intra-block coherence of the dictionary $\Phi$ is defined as $$\mu_I := \max_{1\le i \le r} \|\Phi_i^*\Phi_i- \ident_m\|_2,$$ the inter-block coherence of the dictionary $\Phi$ is defined as\footnote{See \cite{EKB10} for a related measure of block coherence of a dictionary that is given by $\mu_B/m$.} $$\mu_B := \max_{1\le i\ne j \le r} \|\Phi_i^*\Phi_j\|_2,$$ and the quadratic-mean block coherence of the dictionary $\Phi$ is defined as $$\overline{\mu}_B := \max_{1 \le j \le r} \sqrt{\frac{1}{r-1}\sum_{i=1,i\ne j}^r \|\Phi_i^*\Phi_j\|_2^2}.$$
\end{defi}
Note that $\mu_I$ measures the deviation of individual blocks $\{\Phi_i\}$ from being orthonormal and is identically equal to zero for the case of orthonormal blocks. In contrast, $\mu_B$ measures the similarity between different blocks and cannot be zero in the $n$ smaller than $p$ setting. Finally, $\overline{\mu}_B$ measures the \emph{average} similarity between different blocks and is a generalization of the \emph{mean square coherence} defined in \cite{Barg.etal.ITIT2015} to the block setting.

In addition to these three measures, the main result also relies on a condition that we term the \emph{block incoherence condition} (BIC).
\begin{defi}[Block Incoherence Condition]
We say that $\Phi$ satisfies the \emph{block incoherence condition} (BIC) with parameters $(c_1, c_2)$ if $\mu_I \leq c_1$ and $\mu_B \leq c_2/\log{p}$ for some positive numerical constants $c_1$ and $c_2$.
\end{defi}
Informally, the BIC dictates that individual blocks of $\Phi$ do not diverge from being orthonormal in an unbounded fashion and the \emph{worst-case} dissimilarity between different blocks scales as $O(1/\log{p})$. It is worth making two important observations here. First, the BIC does not require the intra- and inter-block coherences of $\Phi$ to scale with the number of blocks, $k$, in the block subdictionary $X$. Second, in contrast to the inter-block coherence $\mu_B$, the quadratic-mean block coherence $\overline{\mu}_B$ of the dictionary does have an effect on the number of blocks $k$ in $X$; hence, it is not included in the BIC. We conclude by noting that the most desirable aspect of the BIC is that it can be verified in polynomial time. We are now ready to state our first result.
\begin{thm}\label{thm:rand_cond}
Suppose that the $n \times p$ dictionary $\Phi = [\Phi_1~\Phi_2~\ldots~\Phi_r]$ satisfies the BIC with parameters $(c_1, c_2)$. Let $\cl{S}$ be a $k$-subset drawn uniformly at random from all ${r \choose k}$ possible $k$-subsets of $\{1,\dots,r\}$ with $k \ge \log p$. Then, as long as $k \leq \min\left\{c_0/(\overline{\mu}_B^2 \log{p}), c_0^\prime r/(\|\Phi\|^2_2 \log{p})\right\}$ for some positive numerical constants $(c_0, c_0^\prime)$ that depend only on $(c_1, c_2)$, the singular values of the block subdictionary $X = [\Phi_i : i \in \cl{S}]$ satisfy $\sigma_i(X) \in [\sqrt{1/2}, \sqrt{3/2}\,]$, $i=1,\dots,km$, with probability with respect to the random choice of the subset $\cl{S}$ of at least $1 - 2p^{-4\log2}$.
\end{thm}
\begin{rem}
The interval $[\sqrt{1/2}, \sqrt{3/2}\,]$ in Theorem~\ref{thm:rand_cond} is somewhat arbitrary. In general, it can be replaced with $[\,\sqrt{1-\epsilon}, \sqrt{1+\epsilon}\,]$ for any $\epsilon \in (0,1)$, resulting in the probability of success either increasing ($\epsilon > 1/2$) or decreasing ($\epsilon < 1/2$).
\end{rem}

In words, Theorem~\ref{thm:rand_cond} states that if a dictionary satisfies the BIC then most of its block subdictionaries of dimensions $n \times km$ act as isometries on $\R^{km}$ for $k = O\left(\min\left\{1/(\overline{\mu}_B^2 \log{p}),r/(\|\Phi\|^2_2 \log{p})\right\}\right)$. In order to better understand the two terms in the bound for $k$, first notice that $\|\Phi\|^2_2 \geq p/n$ for the case of a normalized dictionary \cite{Tropp.ACHA2008}, implying $r/(\|\Phi\|^2_2 \log{p}) = O(n/(m\log{p}))$. More importantly, the equality $\|\Phi\|^2_2 = p/n$ is achievable by dictionaries with orthogonal rows (also referred to as tight frames \cite{Bajwa.Pezeshki.FF2012}). Further, there also exist block dictionaries for which $\overline{\mu}_B^2 = O(m/n)$ \cite{Calderbank.etal.ACHA2015}, implying $1/(\overline{\mu}_B^2 \log{p}) = O(n/(m\log{p}))$. Together, these two facts imply that Theorem~\ref{thm:rand_cond} allows optimal scaling of the dimensions of well-conditioned block subdictionaries for certain dictionaries, i.e., it does not suffer from the square-root bottleneck. Perhaps the most surprising aspect of this theorem, which sets it apart from other works on inference under linear models in block settings \cite{TGS06,T06,GRSV08,EKB10,ER10,Lounici.etal.AS2011,Elhamifar.Vidal.ITSP2012}, is the assertion it makes about the effects of the intra- and inter-block coherences of $\Phi$ on the conditioning of random block subdictionaries. Roughly, Theorem~\ref{thm:rand_cond} suggests that as soon as the BIC is satisfied, both $\mu_I$ and $\mu_B$ stop playing a role in determining the order-wise dimensions of the subdictionaries that are well conditioned; rather, it is the spectral norm of the dictionary $\|\Phi\|_2$ and the quadratic-mean block coherence of the dictionary $\overline{\mu}_B$ that play primary roles in this regard. Such an assertion of course needs to be carefully examined, given that Theorem~\ref{thm:rand_cond} is only concerned with sufficient conditions. Nevertheless, carefully planned numerical experiments carried out in the context of block-sparse recovery (cf.~Section~\ref{sec:grouprecovery}) and block-sparse regression (cf.~Section~\ref{sec:grouplasso}) are consistent with this assertion.

\subsection{Proof of Theorem~\ref{thm:rand_cond}}
The proof of Theorem~\ref{thm:rand_cond} leverages the analytical tools employed by Tropp in \cite{TroppCRAS} for conditioning of canonical (i.e., non-block) random subdictionaries, coupled with a \emph{Poissonization} argument that is now standard in the literature (see, e.g, \cite{CandesPlan}). To proceed, we define $r$ independent and identically distributed (i.i.d.) Bernoulli random variables $\zeta_1,\ldots,\zeta_r$ with parameter $\delta := k/r$ (i.e., $\prob(\zeta_i = 1) = \delta$) and a random set $\cl{S}^\prime := \{i : \zeta_i = 1\}$. Next, we define a random block subdictionary $X^\prime := [\Phi_i : i \in \cl{S}^\prime]$ and use $G := \Phi^*\Phi - \ident$ and $F := {X^\prime}^*X^\prime - \ident$ to denote the hollow Gram matrix of $\Phi$ and the hollow Gram matrix of $X^\prime$, respectively. Finally, define $\Sigma := \textrm{diag}(\zeta_1,\ldots,\zeta_r)$ to be a random diagonal matrix, $R := \Sigma \otimes \ident_m$ to be a block masking matrix, and notice from definition of the spectral norm that $\|F\|_2 = \|RGR\|_2$. Using this notation, we can show that the $L_q$ norm of the random variable $\|RGR\|_2$ for $q = 4\log{p}$ is controlled by $\mu_I$, $\mu_B$, $\overline{\mu}_B$, and $\|\Phi\|_2$.
\begin{lem}
\label{lem:cond_block} For $\delta = k/r$ and $q=4\log{p}$,  the $L_q$ norm of the random variable $\|RGR\|_2 = \|F\|_2$ can be bounded as $$\Eq \|RGR\|_2 \le 48\mu_B\log p +6\overline{\mu}_B\sqrt{2(r-1)\delta \log p} + 6\sqrt{2\delta\log{p}(1+\mu_I)}\|\Phi\|_2 + 2\delta\|\Phi\|_2^2+3\mu_I.$$
\end{lem}

The proof of Lemma~\ref{lem:cond_block}, which is fundamental to the proof of Theorem~\ref{thm:rand_cond} and comprises novel generalizations of some of the results in \cite{TroppCRAS,TroppRPP,Rudelson99,RudelsonVershynin} to the block setting of this paper, is provided in Appendix~\ref{sec:mainproof}. We are now ready to provide a proof of Theorem~\ref{thm:rand_cond}.
\begin{IEEEproof}[Proof of Theorem~\ref{thm:rand_cond}]
Define $Z := \|X^*X-\ident\|_2$ and notice that $\sigma_i(X) \in [\sqrt{1/2}, \sqrt{3/2}\,]$, $i=1,\dots,km$, if and only if $Z \leq 1/2$. Instead of studying $Z$ directly, however, we first study the related random variable $Z^\prime := \|{X^\prime}^*X^\prime-\ident\|_2 = \|F\|_2$, where $X^\prime$ is the random subdictionary defined in relation to Lemma~\ref{lem:cond_block}. It then follows from the Markov inequality and Lemma~\ref{lem:cond_block} that
\begin{align}
\prob(Z' > 1/2) &\le (1/2)^{-q}(\Eq[Z'])^q\nonumber\\
&\le 2^q \left(48\mu_B\log p +6\overline{\mu}_B\sqrt{2(r-1)\delta \log p} + 6\sqrt{2\delta\log{p}(1+\mu_I)}\|\Phi\|_2 + 2\delta\|\Phi\|_2^2+3\mu_I\right)^q,
 \label{eq:zprime_prf}
\end{align}
where $q := 4\log{p}$. Next, our goal is to show that for all $t > 0$,
\begin{align}
\prob(Z > t) \le 2 \prob(Z' > t){,}
\label{eq:poissonization}
\end{align}
using the Poissonization argument from~\cite{CandesPlan}. Toward this end, we explicitly write $X' = \Phi_{\cl{S}'}$ and note that
\begin{align}
\nonumber
\prob(\|\Phi_{\cl{S}'}^*\Phi_{\cl{S}'}-\ident\|_2 > t) &= \sum_{{{\ell}}=0}^r \prob\left(\|\Phi_{\cl{S}'}^*\Phi_{\cl{S}'}-\ident\|_2 > t \big\vert\, |\cl{S}'| = {{\ell}}\right)\prob(|\cl{S}'| = {{\ell}}) \\
\nonumber
&\ge \sum_{{{\ell}}=k}^r \prob\left(\|\Phi_{\cl{S}'}^*\Phi_{\cl{S}'}-\ident\|_2 > t \big\vert\, |\cl{S}'| = {{\ell}}\right)\prob(|\cl{S}'| = {{\ell}}) \\
&= \sum_{{{\ell}}=k}^r \prob\left(\|\Phi_{\cl{S}_{{\ell}}}^*\Phi_{\cl{S}_{{\ell}}}-\ident\|_2 > t \right)\prob(|\cl{S}'| = {{\ell}}),
\end{align}
where $\cl{S}_{{\ell}}$ is a subset drawn uniformly at random from all ${r \choose \ell}$ possible $\ell$-subsets of $\{1,\dots,r\}$. We now make two observations. First, $|\cl{S}'|$ is a binomial random variable with parameters $(r, k/r)$ and therefore $\prob(|\cl{S}'| \ge k) \ge 1/2$ due to $k$ being the median of $|\cl{S}'|$. Second, since each (random) submatrix $\Phi_{\cl{S}_{\ell'}}^*\Phi_{\cl{S}_{\ell'}}-\ident$ for a given value of $\ell' \le \ell$ is a submatrix of some (random) $\Phi_{\cl{S}_{\ell}}^*\Phi_{\cl{S}_{\ell}}-\ident$ and the spectral norm of a matrix is lower bounded by that of its submatrices, we have that $\prob(\|\Phi_{\cl{S}_{\ell}}^*\Phi_{\cl{S}_{\ell}}-\ident\|_2 > t)$ is a nondecreasing function of ${\ell}$. Therefore we can write
\begin{align}
\prob(\|\Phi_{\cl{S}'}^*\Phi_{\cl{S}'}-\ident\|_2 > t) &\ge \prob\left(\|\Phi_{\cl{S}_k}^*\Phi_{\cl{S}_k}-\ident\|_2 > t \right)\sum_{{\ell}=k}^r \prob(|\cl{S}'| = {\ell}) \nonumber \\
&\ge \prob\left(\|\Phi_{\cl{S}_k}^*\Phi_{\cl{S}_k}-\ident\|_2 > t \right)\prob(|\cl{S}'| \ge k) \nonumber \\
&\ge \frac{1}{2}\prob\left(\|\Phi_{\cl{S}_k}^*\Phi_{\cl{S}_k}-\ident\|_2 > t \right) \nonumber \\
&= \frac{1}{2}\prob\left(\|\Phi_\cl{S}^*\Phi_\cl{S}-\ident\|_2 > t \right), \label{eq:poisscp_prf}
\end{align}
where the last equality follows since $\cl{S}_k$ and $\cl{S}$ have the same probability distribution. By combining (\ref{eq:zprime_prf}) and (\ref{eq:poisscp_prf}), we therefore obtain
\begin{align}
\prob(Z > 1/2)&\le 2^{q+1} \left(48\mu_B\log p +6\overline{\mu}_B\sqrt{2(r-1)\delta \log p} + 6\sqrt{2\delta\log{p}(1+\mu_I)}\|\Phi\|_2 + 2\delta\|\Phi\|_2^2+3\mu_I\right)^q.
\end{align}
Finally, the expression inside parentheses in the above equation can be bounded by $1/4$ for small-enough constants $c_0, c_0^\prime, c_1$, and $c_2$, resulting in
$\prob(Z > 1/2) \le 2(1/2)^{4\log{p}} = 2p^{-4\log 2}$.
\end{IEEEproof}

\subsection{Discussion}
\label{sec:mmv}
Among existing works focusing on the conditioning of random (non-block) subdictionaries \cite{Candes04C,CandesQRUP,Tropp.ACHA2008,TroppCRAS,Kuppinger.etal.ITIT2012,Gurevich.Hadani.unpublished2009}, \cite{TroppCRAS} and \cite{Kuppinger.etal.ITIT2012} are the ones with the most general and strongest results. Specifically, \cite{Candes04C,CandesQRUP} deal with the case of the dictionary $\Phi$ being a concatenation of two orthonormal bases, \cite{Gurevich.Hadani.unpublished2009} studies the case of $\Phi$ being a disjoint union of orthonormal bases, and \cite{Tropp.ACHA2008} requires $k$ to be inversely proportional to the column-wise coherence $\mu^2(\Phi)$. The results in \cite{TroppCRAS} and \cite{Kuppinger.etal.ITIT2012} are related to each other in the sense that \cite{Kuppinger.etal.ITIT2012} extends \cite{TroppCRAS} to the case when the subdictionaries of $\Phi$ are not necessarily selected uniformly at random. The proof technique employed in this paper for conditioning of random block subdictionaries is inspired by \cite{TroppCRAS} and is rather tight in the sense that in the case of $m=1$ (and a unit-norm dictionary $\Phi$), we have $\mu_I = 0$, $\overline{\mu}_B \le \frac{\|\Phi\|_2}{\sqrt{p}}$, and therefore Lemma~\ref{lem:cond_block} reduces to \cite[Corollary~5.2]{TroppCRAS}. While we believe our result can be extended to the case when the random block subdictionaries of $\Phi$ are selected with a more ``structured randomness'' by leveraging the insights offered by \cite{Kuppinger.etal.ITIT2012}, we leave this for future work.

It is also instructive to note that while Theorem~\ref{thm:rand_cond} is the most general incarnation of results concerning conditioning of random block subdictionaries, it is rather straightforward to specialize this result for conditioning of random block subdictionaries of \emph{structured} dictionaries. Next, we specialize Theorem~\ref{thm:rand_cond} to one such structure that corresponds to $\Phi$ being a Kronecker product of an arbitrary unit-norm dictionary and a dictionary with orthonormal columns. Such \emph{Kronecker-structured} dictionaries arise in many contexts~\cite{Sun,Duarte.Eldar.ITSP2011} and have a special connection to the literature on multiple measurement vectors (MMV) \cite{Cotter_Sparse,TGS06,T06,GRSV08,FR08,EM09,ER10,DE10,Lee.etal.ITIT2012,Mishali,Kim,DWBSB13} and multivariate linear regression \cite{Turlach2005,ArgyriouNIPS2006,Obozinski2010,ObozinskiUnion} problems. The following section therefore will also help understand our work in the context of these two research areas.

\subsubsection{Random block subdictionaries of Kronecker-structured dictionaries with application to multiple measurement vectors problem}
Consider an arbitrary, unit-norm-column $n_1 \times r$ dictionary $P$ and an $n_2 \times m$ dictionary $Q$ with orthonormal columns (i.e., $Q^*Q = \ident_m$), where $n_1 > r$, $n_2 \leq m$, and $n_1n_2 = n$. Then a corollary of Theorem~\ref{thm:rand_cond} is that conditioning of random block subdictionaries of the Kronecker-structured dictionary $\Phi = P \otimes Q$ is simply a function of the coherence, $\mu(P) = \max_{i,j:i\not=j} |\langle P_i,P_j \rangle|$, and spectral norm, $\|P\|_2$, of $P$, where $P_i$ denotes the $i^{th}$ column of $P$. Formally, this corollary has the following statement.
\begin{cor}\label{cor:kron_rand_cond}
Suppose that the $n \times p$ dictionary $\Phi = P \otimes Q$ with $Q^*Q = \ident_m$ and $\mu(P) \leq c_2/\log{p}$ for a positive numerical constant $c_2$. Let $\cl{S}$ be a $k$-subset drawn uniformly at random from all ${r \choose k}$ possible $k$-subsets of $\{1,\dots,r\}$. Then, as long as $k \leq c_0 r/(\|P\|^2_2 \log{p})$ for some positive numerical constant $c_0 := c_0(c_2)$, the singular values of the block subdictionary $X = [\Phi_i = P_i \otimes Q: i \in \cl{S}]$ satisfy $\sigma_i(X) \in [\sqrt{1/2}, \sqrt{3/2}]$, $i=1,\dots,km$, with probability at least $1 - 2p^{-4\log2}$. Here, the probability is with respect to the random choice of the subset $\cl{S}$.
\end{cor}

Corollary~\ref{cor:kron_rand_cond} is a simple consequence of properties of Kronecker product. In terms of the spectral norm of $\Phi$, we have $\|P \otimes Q\|_2 = \|P\|_2\|Q\|_2 = \|P\|_2$. In terms of the intra- and inter-block coherences, we note that
\begin{align}
\Phi_i^*\Phi_j = (P_i \otimes Q)^*(P_j \otimes Q) = (P_i^* \otimes Q^*)(P_j \otimes Q) = (P_i^*P_j) \otimes (Q^*Q)=\langle P_i,P_j\rangle \ident_m,
\end{align}
which trivially leads to
\begin{align}
\mu_I &= \max_{1\le i \le r} \|\Phi_i^*\Phi_i- \ident_m\|_2 = \max_{1\le i \le r}\|(\langle P_i,P_i\rangle-1)\ident_m\|_2 =0, \quad \text{and}\label{eq:kron2}\\
\mu_B &= \max_{1\le i\ne j \le r} \|\Phi_i^*\Phi_j\|_2 = \max_{1\le i\ne j \le r}|\langle P_i,P_j\rangle| = \mu(P),\label{eq:kron1}\\
\overline{\mu}_B &= \max_{1 \le j \le r} \sqrt{\frac{1}{r-1}\sum_{i=1,i\ne j}^r \|\Phi_i^*\Phi_j\|_2^2} = \max_{1 \le j \le r} \sqrt{\sum_{i=1,i\ne j}^r \frac{|\langle P_i,P_j\rangle|^2}{r-1}} \leq \max_{1 \le j \le r} \frac{\|P^*P_j\|_2}{\sqrt{r-1}} \le \frac{\|P\|_2}{\sqrt{r-1}}.
\end{align}
Note that Corollary~\ref{cor:kron_rand_cond} is not the tightest possible result for Kronecker-structured dictionaries since Theorem~\ref{thm:rand_cond} does not exploit any dictionary structure. In particular, one can obtain a variant of Corollary~\ref{cor:kron_rand_cond} in which the $\log{p}$ terms are replaced with the $\log{r}$ terms by explicitly accounting for the Kronecker structure in the proof of Lemma~\ref{lem:cond_block}.

We conclude this section by connecting Corollary~\ref{cor:kron_rand_cond} to the MMV/multivariate linear regression problem, which will help clarify the similarities and differences between our work and MMV-related works.\footnote{In the operational sense, MMV and multivariate linear regression are two distinct inference problems. For ease of exposition, however, we use the term MMV in here to refer to both problems.} The inference problems studied under the MMV setting are essentially special cases of inference in the block-sparse setting studied here. In the MMV setting, it is assumed there are $m$ parameter vectors, $b_1,\ldots,b_m$, collected as $m$ columns of an $r \times m$ matrix $B$. In addition, each $b_i$ is observed using the same $n \times r$ dictionary $A$, $y_i := A b_i$, and the observation vectors $y_1,\ldots,y_m$ are collected as $m$ columns of an $n \times m$ matrix $Y := A B$. A typical assumption in this MMV setting states that the $m$ different parameter vectors share locations of their $k \ll r$ nonzero entries, resulting in $B$ having no more than $k$ nonzero \emph{rows}. It is however easy to see that if we define $y := \text{vec}(Y^T)$ and $b := \text{vec}(B^T)$ then $y = (A \otimes \ident_m)b$, where the vector $b$ exhibits block sparsity. In other words, inference in the MMV setting requires understanding the conditioning of block subdictionaries of $A \otimes \ident_m$. In this case, we already know from Corollary~\ref{cor:kron_rand_cond} that the conditioning of random block subdictionaries of $A \otimes \ident_m$ is simply a function of the coherence and spectral norm of $A$. Interestingly, while there exists a significant body of literature in the MMV setting \cite{Cotter_Sparse,TGS06,T06,GRSV08,FR08,EM09,ER10,DE10,Lee.etal.ITIT2012,Mishali,Kim,Turlach2005,ArgyriouNIPS2006,Obozinski2010,ObozinskiUnion,DWBSB13}, most of these works do not provide near-optimal, verifiable conditions for guaranteeing success of MMV-based inference problems. The most notable exception to this is the recent work \cite{ER10}, which studies the problem of noiseless recovery in the MMV setting. Nonetheless, our block-sparsity results (including the forthcoming noiseless recovery results) are much more general than the ones in \cite{ER10} because of the MMV setting being just a special case of Corollary~\ref{cor:kron_rand_cond} in the block-sparse setting.

\section{Application: Recovery of Block-Sparse Signals from Noiseless Measurements}
\label{sec:grouprecovery}
We now shift our focus to the applicability of Theorem~\ref{thm:rand_cond} in the context of inference problems. We first begin with the problem of recovery of $\beta$ from $y=\Phi\beta$ when the signal $\beta$ is block sparse. Block sparsity is one of the most popular structures used in sparse signal recovery problems. It is also intrinsically linked with the multiple measurement vectors (MMV) problem described in Section~\ref{sec:mmv}, as there is an equivalent block-sparse formulation for each MMV problem. Block sparsity arises in many applications, including union-of-subspaces models~\cite{EM09,Blumensath.Davies.ITIT2009}, multiband communications~\cite{MEDS09,MishaliEldar}, array processing~\cite{MCW05,GR97}, and multi-view medical imaging~\cite{PLM97,GR97,OuEEG}.

Because of the relevance of block sparsity in these and other applications, significant efforts have been made toward development of block-sparse signal recovery methods/algorithms and matching guarantees on the number of measurements required for successful recovery \cite{Stojnic.etal.ITSP2009,Stojnic2,EM09,T06,TGS06,GRSV08,Elhamifar.Vidal.ITSP2012,fusion,FangLi,EKB10,BenHaim,Cotter_Sparse,Lee.etal.ITIT2012,Kim,DE10,EM09,BCDH10,ER10,Mishali,DWBSB13}. However, the results reported in some of these works are only applicable in the case of randomized dictionary constructions~\cite{Stojnic.etal.ITSP2009,Stojnic2,EM09,BCDH10,DWBSB13}, while those reported in other works rely on dictionary conditions that either cannot be explicitly verified in polynomial time \cite{Cotter_Sparse,TGS06,GRSV08,Elhamifar.Vidal.ITSP2012,Lee.etal.ITIT2012,Kim,DE10,EM09,fusion,BCDH10,Mishali} or result in a suboptimal scaling of the number of measurements due to their focus on the worst-case performance \cite{T06,TGS06,GRSV08,Elhamifar.Vidal.ITSP2012,fusion,FangLi,EKB10,BenHaim}.

To the best of our knowledge, the only work that does not have the aforementioned limitations is \cite{ER10}. Nonetheless, the focus in \cite{ER10} is only on the restrictive MMV problem, rather than the general block-sparse signal recovery problem. In addition, the analytical guarantees provided in \cite{ER10} rely on the nonzero entries of $\beta$ following either Gaussian or spherical distributions. In contrast, we make use of the main result of Section~\ref{sec:result} in the following to state a result for average-case recovery of block-sparse signals that suffers from none of these and earlier limitations. Our result depends primarily on the spectral norm and the quadratic-mean block coherence of $\Phi$, while it has a mild dependence on the intra- and inter-block coherences; and all four of these quantities can be explicitly computed in polynomial time. It further requires only weak assumptions on the distribution of the nonzero entries of $\beta$. Equally important, the forthcoming result does not suffer from the so-called ``square-root bottleneck'' \cite{Kuppinger.etal.ITIT2012}; specifically, for the case when $m = O(1)$, the stated result allows almost linear scaling of the block-sparsity level $k$ as a function of the number of measurement $n$ for dictionaries $\Phi$ with sufficiently small spectral norms (e.g., tight frames) and quadratic-mean block coherences.

\subsection{Recovery of Block-Sparse Signals: Problem Formulation}\label{ssec:stat_model}
Our exposition throughout the rest of this section will be based upon the following formulation. We are interested in recovering a block-sparse signal $\beta \in \R^p$ from noiseless measurements $y = \Phi \beta$, where the dictionary $\Phi$ denotes an $n \times p$ observation matrix with $n \ll p$ and $y \in \R^n$ denotes the observation vector. We assume $\beta$ comprises a total of $r$ blocks, each of size $m$ (yielding $p = rm$), and represent it without loss of generality as $\beta = [\beta_1^*~\beta_2^*~\ldots~\beta_r^*]^*$ with each block $\beta_i \in \R^m$. In order to make this problem well posed, we require that $\beta$ is $k$-block sparse with $\#\{i:\beta_i \not= {\bf 0}\} = k \ll r$. Finally, we impose a mild statistical prior on $\beta$, as described below.
\begin{enumerate}
\item[M1)] The \emph{block support} of $\beta$, $\cl{S} = \{i:\beta_i \not= {\bf 0}\}$, has a uniform distribution over all $k$-subsets of $\{1,\dots,r\}$,
\item[M2)] The entries in $\beta$ have zero median (i.e., the nonzero entries are equally likely to be positive and negative): $\mathbb{E}(\sign(\beta)) = {\bf 0}$, and
\item[M3)] Nonzero blocks of the block-sparse signal $\beta$ have statistically independent ``directions.'' Specifically, we require $\mathbb{P}\left(\bigcap_{i\in \cl{S}}\left(\signbar(\beta_i) \in \mathcal{A}_i\right)\right) = \prod_{i\in \cl{S}}\mathbb{P} \left(\signbar(\beta_i) \in \mathcal{A}_i\right)$, where $\mathcal{A}_i \subset \mathbb{S}^{m-1}$ with $\mathbb{S}^{m-1}$ denoting the unit sphere in $\real^m$.
\end{enumerate}
Note that M2 and M3 are trivially satisfied in the case of the nonzero blocks of $\beta$ drawn independently from either Gaussian or spherical distributions. However, it is easy to convince oneself that many other distributions---including those that are not absolutely continuous---will satisfy these two conditions. Conditions M1--M3 provide a probabilistic characterization of block-sparse $\beta$ that is inspired by Tropp \cite{Tropp.ACHA2008} and Cand\`{e}s and Plan \cite{CandesPlan} in which a related characterization of non block-sparse $\beta$ helped them overcome some analytical hurdles in relation to performance specifications of sparse recovery and regression problems, respectively.

\subsection{Main Result and Discussion}
In this section, we are interested in understanding the average-case performance of the following mixed-norm convex optimization program for recovery of block-sparse signals satisfying M1--M3:
\begin{align}
\betahat = \arg \min_{\bar{\beta} \in \real^p} \|\bar{\beta}\|_{2,1}~\textrm{such that}~y = \Phi \bar{\beta},
\label{eq:l21min}
\end{align}
where the $\ell_{2,1}$ norm of a vector $\beta \in \real^{p}$ containing $r$ blocks of $m$ entries each is defined as $\|\beta\|_{2,1} := \sum_{i=1}^{{r}} \|\beta_i\|_2$. While \eqref{eq:l21min} has been utilized in the past for recovery of block sparse signals (see, e.g., \cite{EM09,Stojnic.etal.ITSP2009,ER10}), an average-case analysis result along the following lines is novel. The following theorem is proven in Appendix~\ref{sec:l12proof}.
\begin{thm}
Suppose that $\beta \in \R^p$ is $k$-block sparse and it is drawn according to the statistical model M1, M2, and M3. Further, assume that $\beta$ is observed according to the linear model $y = \Phi \beta$, where the $n \times p$ matrix $\Phi$ satisfies $\mu_I \leq c_1$ and $\mu_B \leq c_2'/(\sqrt{m} \log p)$ for some positive numerical constants $c_1$ and $c_2'$. Then, as long as $k \leq \min\left\{c_0''/(\overline{\mu}_B^2 m \log{p}), c_0^\prime r/(\|\Phi\|^2_2 \log{p})\right\}$ for some positive numerical constants $(c_0'', c_0')$ that depend only on $(c_1, c_2')$, the minimization (\ref{eq:l21min}) results in $\betahat = \beta$ with probability at least $1-6p^{-4\log 2}$. \label{thm:l12}
\end{thm}

It is worth pointing out here that the conditions on $\mu_I$ and $\mu_B$ stated in Theorem~\ref{thm:l12} imply the BIC. These stronger conditions on the inter- and intra-block coherences will also be used in several of our forthcoming results. Interestingly, Theorem~\ref{thm:l12} specialized to the case of non-block sparse signals (by setting $m = 1$ and $r = p$) gives us an average-case analysis result for recovery of sparse signals that has never been explicitly stated in prior works. The optimization program \eqref{eq:l21min} in this case reduces to the standard \emph{basis pursuit} program \cite{DonohoBP}:
\begin{align}
\betahat = \arg \min_{\bar{\beta} \in \real^{p}} \|\bar{\beta}\|_1~\textrm{such that}~y = \Phi \bar{\beta},
\label{eq:l1}
\end{align}
the BIC reduces to a bound on the coherence of $\Phi$, and Theorem~\ref{thm:l12} reduces to the following corollary.\footnote{We refer the reader to the forthcoming discussion for the difference between the average-case analysis result in \cite{Tropp.ACHA2008} and Corollary~\ref{thm:l1}. While it is possible to leverage the results in \cite{TroppCRAS} to obtain Corollary~\ref{thm:l1}, rather than obtaining Corollary~\ref{thm:l1} from Theorem~\ref{thm:l12} in this paper, such a result does not explicitly exist in prior literature to the best of our knowledge.}
\begin{cor}
Suppose $\beta \in \R^p$ is $k$-sparse, its support (i.e., locations of its nonzero entries) is a $k$-subset drawn uniformly at random from all ${p \choose k}$ possible $k$-subsets of $\{1,\dots,p\}$, its nonzero entries are drawn from a multivariate distribution with zero median (i.e., the nonzero entries are equally likely to be positive and negative), and the signs of the nonzero entries are independent. Then, as long as $k \le c_0' p/\|\Phi\|_2^2 \log p$ and $\mu(\Phi) \le c_1'/\log p$ for some positive numerical constants $c_0' := c_0'(c_1')$ and $c_1'$, the minimization (\ref{eq:l1}) successfully recovers $\beta$ from $y = \Phi \beta$ with probability at least $1-6p^{-4\log 2}$.\label{thm:l1}
\end{cor}

We now elaborate on the similarities and differences between our (average-case) guarantees for recovery of block-sparse (Theorem~\ref{thm:l12}) and non-block sparse (Corollary~\ref{thm:l1}) signals. In terms of similarities, if we assume $k = O(1/\overline{\mu}_B^2 m \log{p})$ then both results allow for the same scaling of the \emph{total number of nonzero entries} in $\beta$: $km = O(p/\|\Phi\|_2^2 \log p)$ in the case of block-sparse signals and $k = O(p/\|\Phi\|_2^2 \log p)$ in the case of sparse signals. However, while Corollary~\ref{thm:l1} requires that the inner product of any two columns in $\Phi$ be $O(1/\log p)$, Theorem~\ref{thm:l12} allows for less restrictive inner products of columns \emph{within} blocks as long as $\mu_I = O(1)$ and $\mu_B = O(1/\sqrt{m} \log p)$. 
Similarly, while Corollary~\ref{thm:l1} requires that the signs of the nonzero entries in $\beta$ be independent, Theorem~\ref{thm:l12} allows for correlations among the signs of entries \emph{within} nonzero blocks. With the caveat that Theorem~\ref{thm:l12} and Corollary~\ref{thm:l1} only specify sufficient conditions, these two results seem to suggest that explicitly accounting for block structures in sparse signals with $m = O(1)$ allows one to expand the classes of sparse signals $\beta$ \emph{and} dictionaries $\Phi$ under which successful (average-case) recovery can be guaranteed.

Next, we comment on the scaling of the number of nonzero entries in Theorem~\ref{thm:l12}. Assuming the size of the blocks $m = O(1)$ and appropriate conditions on statistical properties of $\beta$ and intra-/inter-block coherences of $\Phi$ are satisfied, Theorem~\ref{thm:l12} allows for the number of nonzero entries to scale like $O(n/\log p)$ for dictionaries $\Phi$ that satisfy $\overline{\mu}_B^2 = O(m/n)$ \cite{Calderbank.etal.ACHA2015} and that are ``approximately'' tight frames \cite{Bajwa.Pezeshki.FF2012}, $\|\Phi\|_2^2 \approx p/n$. This suggests a near-optimal nature of Theorem~\ref{thm:l12} (modulo perhaps $\log$ factors) for sufficiently well-behaved dictionaries as one cannot expect better than linear scaling of the number of nonzero entries as a function of the number of observations. In particular, existing literature on frame theory \cite{BCM12} can be leveraged to specialize the results of Theorem~\ref{thm:l12} and Corollary~\ref{thm:l1} for oft-used classes of random dictionaries (e.g., Gaussian, random partial Fourier) and to establish that in such cases the scaling of our guarantees matches that obtained using nonverifiable conditions such as the restricted isometry property \cite{Rudelson.Vershynin.CPAM2008,CandesRIP}.

Additionally, we note that when Corollary~\ref{thm:l1} is specialized to the case of (approximately) tight frames we obtain average-case guarantees somewhat similar to the ones reported in \cite[Theorem 14]{Tropp.ACHA2008}. The main difference between the two results is the role that the coherence $\mu(\Phi)$ plays in the guarantees. In \cite[Theorem 14]{Tropp.ACHA2008}, the maximum allowable sparsity $k$ is required to be inversely proportional to $\mu^2(\Phi)$. In contrast, we assert that the maximum allowable sparsity scaling is not fundamentally determined by the (column-wise) coherence $\mu(\Phi)$. Numerical experiments reported in the following section suggest that this is indeed the case.

\begin{rem}
Despite the order-wise tightness of Theorem~\ref{thm:l12} for certain dictionaries with $m = O(1)$, note that the bound on $k$ in it is only a sufficient condition. In particular, we do not have a converse that shows the impossibility of recovering $\beta$ from $y$ when this bound is violated. To the best of our knowledge, however, no such converses exist even in the non block-sparse setting for arbitrary dictionaries and/or non-asymptotic analysis (cf.~\cite{DonohoMar2005,counting}).
\end{rem}

\subsection{Numerical Experiments}\label{ssec:recovery_num_exp}
One of the fundamental takeaways of this section is that the spectral norm of the dictionary, rather than the column-wise coherence of the dictionary, determines the maximum allowable sparsity in (block)-sparse signal recovery problems. In order to experimentally verify this insight, we performed a set of block-sparse signal recovery experiments with carefully designed dictionaries having varying spectral norms and coherence values. Throughout our experiments, we set the signal length to $p=5000$, the block size and the number of blocks to $m = 10$ and $r = 500$, respectively, and the number of observations to $n = 858$ (computed from the bound in~\cite{RRN12} for $k=20$ nonzero blocks). In order to design our dictionaries, we first used Matlab's random number generator to obtain 2000 matrices with unit-norm columns. Next, we manipulated the singular values of each of these matrices to increase their spectral norms by a set of integer multipliers $\mathcal{T}$. Finally, for each of the $2000\cdot |\mathcal{T}|$ resulting matrices, we normalized their columns to obtain our dictionaries and recorded their spectral norms $\|\Phi\|_2$, (column-wise) coherences $\mu(\Phi)$, inter-block coherences $\mu_B(\Phi)$, intra-block coherences $\mu_I(\Phi)$, and quadratic-mean block coherences $\overline{\mu}_B(\Phi)$.

We evaluate the block-sparse signal recovery performance of each resulting dictionary $\Phi$ using Monte Carlo trials, corresponding to the generation of 1000 block-sparse signals with $k$ nonzero blocks. Each signal has block support selected uniformly at random according to M1 and nonzero entries drawn independently from the standard Gaussian distribution $\mathcal{N}(0,\ident)$. We then obtain the observations $y = \Phi \beta$ using the dictionary $\Phi$ under study for each one of these signals and perform recovery using the minimization (\ref{eq:l21min}).\footnote{We used the {\tt SPGL1} Matlab package \cite{spgl1:2007} in all simulations in this section.} We define successful recovery to be the case when the block support of $\betahat$ matches the block support of $\beta$ and the submatrix of $\Phi$ with columns corresponding to the block support of $\beta$ has full rank.

\begin{figure}
\begin{center}
\includegraphics[width=3in]{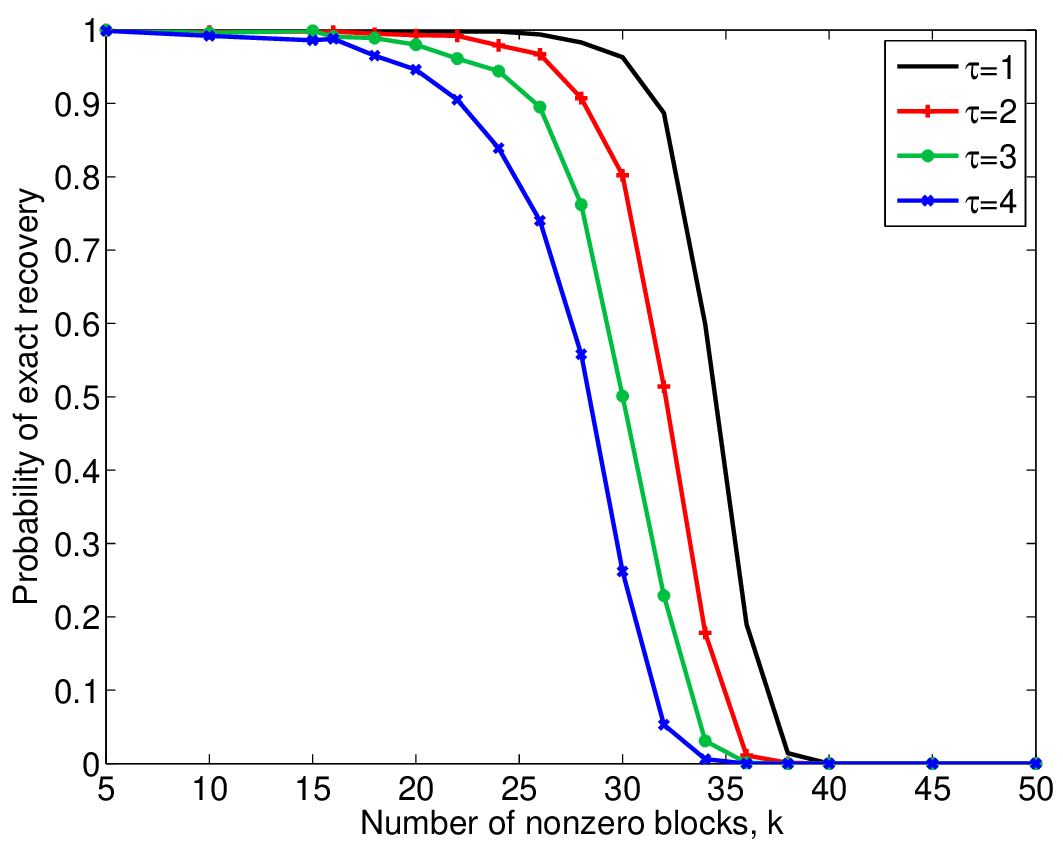}
\end{center}
\caption{Performances of dictionaries $\Phi$ with varying spectral norms and roughly equal coherences in block-sparse signal recovery as a function of the number of nonzero blocks $k$; $\tau \in \mathcal{T}$ denotes the value of the spectral norm multiplier used to generate the dictionary.}
\label{fig:same_mu}
\end{figure}

Figure~\ref{fig:same_mu} shows the performances of dictionaries $\Phi$ of increasing spectral norms ($\mathcal{T} = \{1,2,3,4\}$), where we choose the dictionary (among the 2000 available options) whose coherence value is closest to 0.2. The spectral norms, coherences, inter-block coherences, intra-block coherences, and quadratic-mean block coherences for these chosen dictionaries are collected in Table~\ref{table:same_mu}. The performance is shown as a function of the number of nonzero blocks $k$ in the signal. The figure shows a consistent improvement in the values of $k$ for which successful recovery is achieved as the spectral norm of the dictionary decays, even though $\mu(\Phi)$ does not significantly change among the dictionaries.

\begin{table}[t]
\begin{center}
\begin{tabular}{|l|cccc|}
\hline
$\tau$ & 1 & 2 & 3 & 4 \\
\hline
$\|\Phi_\tau\|_2$ & 3.3963 & 6.7503 & 10.0547 & 13.2034  \\
\hline
$\mu(\Phi_\tau)$ & 0.1992 &   0.2026 &   0.2000  &  0.2207  \\
\hline
$\mu_B(\Phi_\tau)$ & 0.2973 &   0.3431 &   0.5573 &   0.8490 \\
\hline
$\overline{\mu}_B(\Phi_\tau)$ & 0.1965 & 0.2177 & 0.3562 & 0.5624 \\
\hline
$\mu_I(\Phi_\tau)$ & 0.2709 &   0.3131 &   0.4589  &  0.7541 \\
\hline
\end{tabular}
\end{center}
\caption{Spectral norms and coherences for the dictionaries used in the experiments of Figure~\ref{fig:same_mu}.} \label{table:same_mu}
\end{table}

\begin{figure}
\begin{center}
\includegraphics[width=3in]{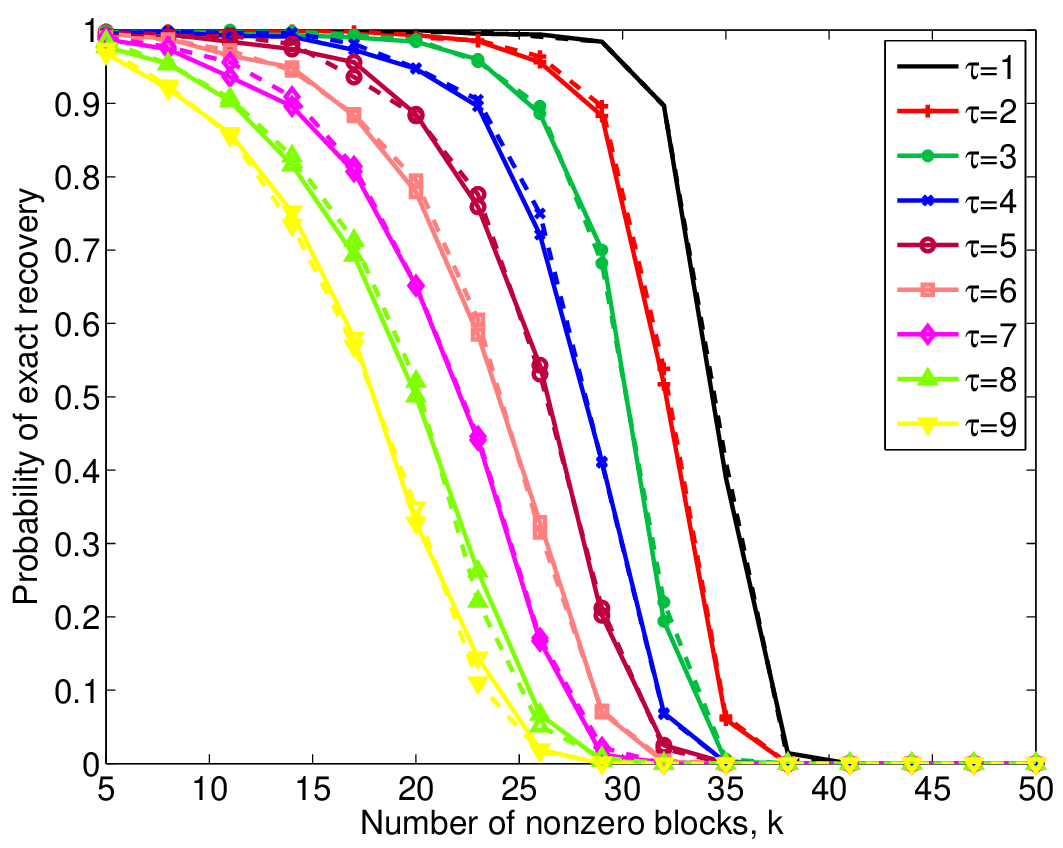}
\end{center}
\caption{Performances of dictionaries $\Phi$ with varying spectral norms and extremal coherence values (cf. Table~\ref{table:extreme_mu}) in block-sparse signal recovery as a function of the number of nonzero blocks $k$; $\tau \in \mathcal{T}$ denotes the value of the spectral norm multiplier used. Solid lines correspond to dictionaries with minimum coherence, while dashed lines correspond to dictionaries with maximum coherence.}
\label{fig:extreme_mu}
\end{figure}

To further emphasize strong dependence of sparse-signal recovery on spectral norm and weak dependence on column-wise coherence, Figure~\ref{fig:extreme_mu} shows the performance of dictionaries $\Phi$ with increasing spectral norms ($\mathcal{T} = \{1,\ldots,9\}$), where we choose dictionaries with the largest and smallest coherence values for each $\tau \in \mathcal{T}$ (among the 2000 available options). The spectral norms, coherences, inter-block coherences, intra-block coherences, and quadratic-mean block coherences for these chosen dictionaries are collected in Table~\ref{table:extreme_mu}. The figure shows not only the same consistent improvement as the spectral norm of the dictionary decays, but also that significant changes in the values of the (column-wise) coherence do not significantly affect the recovery performance. This behavior agrees with our expectation from Theorem~\ref{thm:l12} that the role of the column-wise coherence is decoupled from the scaling of the number of nonzero blocks $k$ (equivalently, number of nonzero entries $km$) in the signal.
\begin{table}[t]
\begin{center}
\begin{tabular}{|l|ccccccccc|}
\hline
$\tau$ & 1 & 2 & 3 & 4 & 5 & 6 & 7 & 8 & 9\\
\hline
$\|\Phi_{\tau,{\min}}\|_2$ & 3.4064  &  6.7726  & 10.0536 &  13.2034 &  16.3421 &  19.2980 &  22.1413 &  24.6710 &  27.2951\\
$\|\Phi_{\tau,{\max}}\|_2$ & 3.3963   & 6.7503 &  10.0543 &  13.2250 &  16.2747 &  19.1506 &  21.9975 &  24.7026 &  27.3199\\
\hline
$\mu(\Phi_{\tau,{\min}})$ & 0.1230 &   0.1198 &   0.1500  &  0.2207  &  0.2964  &  0.3760  &  0.4583  &  0.5337  &  0.6000\\
$\mu(\Phi_{\tau,{\max}})$ & 0.1992 &  0.2026  &  0.2698 &  0.3816 &  0.4863 &  0.5778  & 0.6566 &  0.7225  &  0.7758\\
\hline
$\mu_B(\Phi_{\tau,{\min}})$ &
0.2887   & 0.3177 &  0.5357 &  0.8490 &  1.2917 &  1.7372 &  2.2263 & 2.5989 &  3.1204 \\
$\mu_B(\Phi_{\tau,{\max}})$ &
    0.2973  &  0.3431  &  0.6516 &  1.0287  & 1.4419 &  1.6616 &  2.0230 &  2.4479   & 2.8737 \\
\hline
$\overline{\mu}_B(\Phi_{\tau,{\min}})$ & 0.1965 & 0.2187 & 0.3401 & 0.5624 & 0.8451 & 1.1544 & 1.4911 & 1.7809 & 2.1614 \\
$\overline{\mu}_B(\Phi_{\tau,{\max}})$ & 0.1965 & 0.2177 & 0.3843 & 0.6086 & 0.8736 & 1.1405 & 1.4025 & 1.7230 & 2.0824 \\
\hline
$\mu_I(\Phi_{\tau,{\min}})$ & 0.1487 &   0.2002 &   0.3368 &   0.3787 &   0.3472 &   0.4385 &   0.5462 &   1.0551 &   1.3095  \\
$\mu_I(\Phi_{\tau,{\max}})$ & 0.1992 &   0.2026 &   0.2698 &   0.3816 &   0.4863 &   0.5778 &   0.8273 &   1.0415 &   1.2723 \\
\hline
\end{tabular}
\end{center}
\caption{Spectral norms and block coherences for the dictionaries used in the experiments of Figure~\ref{fig:extreme_mu} and Figure~\ref{fig:reg_extreme_mu}.} \label{table:extreme_mu}
\end{table}

\section{Application: Linear Regression of Block-Sparse Vectors}
\label{sec:grouplasso}
In this section, we leverage Theorem~\ref{thm:rand_cond} to obtain average-case results for linear regression of block-sparse vectors, defined as estimating $\Phi \beta$ from $y = \Phi \beta + \text{noise}$ when $\beta$ has a block-sparse structure. In particular, we focus on two popular convex optimization-based methods, the lasso~\cite{LASSO} and the group lasso~\cite{YuanLin}, for characterizing results for linear regression of block-sparse vectors. Our focus on these two methods is due to their widespread adoption by the signal processing and statistics communities. In the signal processing literature, these methods are typically used for efficient sparse approximation{s} of arbitrary signals in overcomplete dictionaries. In the statistics literature, they are mostly used for efficient variable selection and reliable regression under the linear model assumption. Further, ample empirical evidence in both fields suggests that an appropriately regularized group lasso can outperform the lasso whenever there is a natural grouping of the dictionary atoms/regression variables in terms of their contributions to the observations~\cite{YuanLin,BachGroup}. In this section, we analytically characterize the linear regression performances of both the lasso and the group lasso for block-sparse vectors.

Note that analytical characterization of the group lasso using $\ell_1/\ell_2$ regularization for the ``underdetermined'' setting, in which one can have far more regression variables than observations ($n \ll p$), has received attention recently in the statistics literature \cite{BachGroup,MeierGroup,LiuZhang,nardi:ejs08,ZhangGroup,ObozinskiUnion,Kolar,Fang,Lounici.etal.AS2011}. However, prior analytical work on the performance of the group lasso either studies an asymptotic regime~\cite{BachGroup,LiuZhang,MeierGroup,nardi:ejs08}, focuses on random design matrices (i.e., dictionaries)~\cite{BachGroup,MeierGroup,ObozinskiUnion,Kolar}, and/or relies on conditions that are either computationally prohibitive to verify~\cite{LiuZhang,ZhangGroup,nardi:ejs08,Fang} or that do not allow for near-optimal scaling of the number of observations with the number of active blocks of regression variables $k$~\cite{Lounici.etal.AS2011}. In contrast, our analysis for the regression performance of the group lasso using $\ell_1/\ell_2$ regularization in the underdetermined case for block-sparse vectors circumvents these shortcomings of existing works by adopting a probabilistic model, described by the conditions M1--M3 in Section~\ref{ssec:stat_model}, for the blocks of regression coefficients in $\beta$. To the best of our knowledge, the result stated in the sequel concerning the linear regression performances of the group lasso\footnote{We refer to the group lasso using $\ell_1/\ell_2$ regularization as ``group lasso'' throughout the rest of this paper for brevity.} for block-sparse vectors is the first one for block linear regression that is non-asymptotic in nature and applicable to arbitrary design matrices through verifiable conditions, while still allowing for near-optimal scaling of the number of observations with the number of blocks of nonzero regression coefficients for sufficiently well-behaved dictionaries. Our proof techniques are natural extensions of the ones used in \cite{CandesPlan} for the non-block setting and rely on Theorem~\ref{thm:rand_cond} for many of the key steps.

\subsection{Regression of Block-Sparse Vectors: Problem Formulation}
This section concerns regression in the ``underdetermined'' setting for the case when the observations $y \in \R^n$ can be approximately explained by a linear combination of a small number of blocks ($k < n \ll p$) of regression variables (predictors). Mathematically, we have that $y = \Phi \beta + z$, where $\Phi$ denotes the design matrix (dictionary) containing one regression variable per column, $\beta \in \mathbb{R}^p = [\beta_1^*~\beta_2^*~\ldots~\beta_r^*]^*$ denotes the $k$-block sparse vector of regression coefficients corresponding to these variables (i.e., $\#\{i:\beta_i \not= {\bf 0}\} = k \ll r$), and $z \in \R^n$ denotes the modeling error. Here, we assume without loss of generality that $\Phi$ has unit-norm columns, while we assume the modeling error $z$ to be an i.i.d.\ Gaussian vector with variance $\sigma^2$. Finally, in keeping with the earlier discussion, we impose a mild statistical prior on the vector of regression coefficients $\beta$ that is given by the conditions M1, M2, and M3 in Section~\ref{ssec:stat_model}. The fundamental goal in here then is to obtain an estimate $\betahat$ from the observations $y$ such that $\Phi\betahat$ is as close to $\Phi\beta$ as possible, where the closeness is measured in terms of the $\ell_2$ regression error, $\|\Phi\beta - \Phi\betahat\|_2$.

\subsection{Main Results and Discussion}
In this section, we are interested in understanding the average-case regression performance of two methods in the block-sparse setting. The first one of these methods is the lasso \cite{LASSO}, which ignores any grouping of the regression variables and estimates the vector of regression coefficients as
\begin{align}
\label{eq:lasso}
\betahat = \arg \min_{\beta \in \real^p} \frac{1}{2}\|y-\Phi\beta\|_2^2 + 2\lambda\sigma\|\beta\|_1,
\end{align}
where $\lambda > 0$ is a tuning parameter. In terms of a baseline result for the lasso, we can extend the probabilitic model of Cand\`{e}s and Plan \cite{CandesPlan} for non-block linear regression to the block setting and state the following theorem that follows from Theorem~\ref{thm:rand_cond} and Lemma~\ref{lem:cond_block} in this paper and the proof of \cite[Theorem~1.2]{CandesPlan}.
\begin{thm}[\!\!\mbox{\cite[Theorem~1.2]{CandesPlan}}, Theorem~\ref{thm:rand_cond}, and Lemma~\ref{lem:cond_block}]
Suppose that the vector of regression coefficients $\beta \in \R^p$ is $k$-block sparse and that the observation vector can be modeled as $y = \Phi\beta + z$ with the modeling error $z$ being i.i.d.\ Gaussian with variance $\sigma^2$. Further, assume that $\beta$ is drawn according to the statistical model M1 and M2 with the signs of its nonzero entries being i.i.d., and the $n \times p$ matrix $\Phi$ satisfies the BIC with some parameters $(c_1', c_2')$. Then, as long as $k \leq \min\left\{c_0''/(\overline{\mu}_B^2 \log{p}), c_0''' r/(\|\Phi\|^2_2 \log{p})\right\}$ for some positive numerical constants $c_0'' := c_0''(c_1', c_2')$ and $c_0''' := c_0'''(c_1', c_2')$, the lasso estimate {$\betahat$} in \eqref{eq:lasso} computed with $\lambda = \sqrt{2\log p}$ obeys $$\|\Phi\beta-\Phi\betahat\|_2^2 \le C'mk\sigma^2\log p$$ with probability at least $1-O(p^{-1})$, where $C' > 0$ is a constant independent of the problem parameters. \label{thm:lasso_ind.signs}
\end{thm}

The proof of this theorem is omitted here because of its similarity to the proof of the next theorem. While this theorem suggests that the lasso solution in the block setting enjoys many of the optimality properties of the lasso solution in the non-block setting (see, e.g., the discussion in \cite{CandesPlan}), it has one shortcoming. Instead of the less restrictive condition M3, it assumes independence of the signs of the nonzero regression coefficients. In the following, we provide an extension of Theorem~\ref{thm:lasso_ind.signs} to the case of regression of block-sparse vectors with arbitrarily correlated blocks.
The proof of this theorem is given in Appendix~\ref{app:thm:lasso}.
\begin{thm}
Suppose that the vector of regression coefficients $\beta \in \R^p$ is $k$-block sparse and it is drawn according to the statistical model M1, M2, and M3. Further, assume that the observation vector can be modeled as $y = \Phi\beta + z$, where the $n \times p$ matrix $\Phi$ satisfies $\mu_I \leq c_1'$ and $\mu_B \leq c_2''/(\sqrt{m} \log p)$ for some positive numerical constants $c_1'$, $c_2''$, and the modeling error $z$ is i.i.d.\ Gaussian with variance $\sigma^2$. Then, as long as $k \leq \min\left\{c_0''/(\overline{\mu}_B^2 m \log{p}), c_0''' r/(\|\Phi\|^2_2 \log{p})\right\}$ for some positive numerical constants $(c_0'', c_0''')$ that depend only on $(c_1', c_2'')$, the lasso estimate {$\betahat$} in \eqref{eq:lasso} computed with $\lambda = \sqrt{2\log p}$ obeys $$\|\Phi\beta-\Phi\betahat\|_2^2 \le C''mk\sigma^2\log p$$ with probability at least $1-p^{-1}(2\pi \log p)^{-1/2}-8p^{-4\log2}$, where $C'' > 0$ is a constant independent of the problem parameters.
\label{thm:lasso}
\end{thm}

It can be seen from Theorems~\ref{thm:lasso_ind.signs} and \ref{thm:lasso} that both theorems guarantee the same scaling of the regression error, despite the fact that Theorem~\ref{thm:lasso} allows for arbitrary correlations within blocks. Note, however, that the scalings of the maximum number of allowable nonzero blocks and the block coherence in Theorem~\ref{thm:lasso} match the ones in Theorem~\ref{thm:lasso_ind.signs} only for the case of $m=O(1)$. Otherwise, Theorem~\ref{thm:lasso} with correlated blocks will require a more stricter scaling of $\mu_B(\Phi)$, while its scaling of $k$ will be more restrictive if $r/\|\Phi\|^2_2$ dominates $1/(\overline{\mu}_B^2 m)$. It can be seen from the proof of Theorem~\ref{thm:lasso} in Appendix~\ref{app:thm:lasso} that this dependence upon $m$---the size of the blocks---is a direct consequence of allowing for arbitrary correlations within blocks.

Next, we investigate whether it is possible to achieve a similar result using Theorem~\ref{thm:rand_cond} when one explicitly accounts for the block structure of the regression problem. Specifically, the group lasso explicitly accounts for the grouping of the regression variables in its formulation and estimates the vector of regression coefficients as
\begin{align}
\betahat = \arg \min_{\beta \in \real^p} \frac{1}{2}\|y-\Phi\beta\|_2^2 + 2\lambda\sigma\sqrt{m}\|\beta\|_{2,1},
\label{eq:grouplasso}
\end{align}
where $\lambda > 0$ is once again a tuning parameter. The following theorem shows that the group lasso can also achieve the same scaling results as the lasso for block-sparse vectors with arbitrary correlations within blocks (cf.~Theorem~\ref{thm:lasso}). The proof of this theorem is given in Appendix~\ref{sec:grouplassoproof}.
\begin{thm}
Suppose that the vector of regression coefficients $\beta \in \R^p$ is $k$-block sparse and it is drawn according to the statistical model M1, M2, and M3. Further, assume that the observation vector can be modeled as $y = \Phi\beta + z$, where the $n \times p$ matrix $\Phi$ satisfies $\mu_I \leq c_1'$ and $\mu_B \leq c_2''/(\sqrt{m} \log p)$ for some positive numerical constants $c_1'$, $c_2''$, and the modeling error $z$ is i.i.d.\ Gaussian with variance $\sigma^2$. Then, as long as $k \leq \min\left\{c_0''/(\overline{\mu}_B^2 m \log{p}), c_0''' r/(\|\Phi\|^2_2 \log{p})\right\}$ for some positive numerical constants $(c_0'', c_0''')$ that depend only on $(c_1', c_2'')$, the group lasso estimate {$\betahat$} in \eqref{eq:grouplasso} computed with $\lambda = \sqrt{2\log p}$ obeys $$\|\Phi\beta-\Phi\betahat\|_2^2 \le C''mk\sigma^2\log p$$ with probability at least $1-p^{-1}(2\pi \log p)^{-1/2}-8p^{-4\log2}$, where $C'' > 0$ is a constant independent of the problem parameters.
\label{thm:grouplasso}
\end{thm}
\begin{rem}
Note that if one has $m=1$ then block sparsity reduces to the canonical sparsity, the group lasso \eqref{eq:grouplasso} reduces to the lasso \eqref{eq:lasso}, the block coherence $\mu_B(\Phi)$ reduces to the coherence $\mu(\Phi)$, the quadratic-mean block coherence $\overline{\mu}_B \le \|\Phi\|_2$, and Theorem~\ref{thm:grouplasso} essentially reduces to \cite[Theorem~1.2]{CandesPlan}.
\end{rem}

With the caveat that both Theorems~\ref{thm:lasso} and \ref{thm:grouplasso} are concerned with sufficient conditions for average-case regression performance, we now comment on the strengths and weaknesses of these two results. Assuming appropriate conditions are satisfied for the two theorems, we have that both the lasso and the group lasso result in the same scaling of the regression error, $\|\Phi\beta-\Phi\betahat\|_2^2 = O(mk\sigma^2\log p)$, even in the presence of intra-block correlations. This scaling of the regression error is indeed the best that any method can achieve, modulo the logarithmic factor, since we are assuming that the observations are described by a total of $mk$ regression variables. Further, if we assume $k = O(1/\overline{\mu}_B^2 m \log{p})$ then both the lasso and the group lasso allow for near-optimal scaling of the maximum number of regression variables contributing to the observations, $km = O(p/\|\Phi\|_2^2 \log p)$. In fact, similar to the discussion in Section~\ref{sec:grouprecovery}, it is easy to conclude that this scaling of the number of nonzero regression coefficients is near-optimal since it leads to a linear relationship (modulo logarithmic factors) between the number of observations $n$ and the number of active regression variables $km$ for the case of design matrices that are approximately tight frames: $\|\Phi\|_2^2 \approx p/n$. Similar to the case of Theorem~\ref{thm:lasso}, the main difference between Theorems~\ref{thm:lasso_ind.signs} and \ref{thm:grouplasso} is that the scaling of $k$ in Theorem~\ref{thm:grouplasso} will be more restrictive if $r/\|\Phi\|^2_2$ dominates $1/(\overline{\mu}_B^2 m)$. This, and the fact that Theorem~\ref{thm:grouplasso} results in the same guarantees as Theorem~\ref{thm:lasso} even though empirical evidence suggests otherwise, are reasons to believe that better proof techniques can possibly be leveraged to improve upon our stated results for the group lasso. We leave such an investigation for future work.

\begin{rem}[Beyond canonical block sparsity] While the
block-sparse structures of Theorems~\ref{thm:rand_cond}--\ref{thm:grouplasso} can be found in many applications, there exist other applications in which canonical block sparsity does not adequately capture the sparsity structure of $\beta$. Consider, for instance, wavelet expansions of piecewise smooth signals. Nonzero wavelet coefficients in these cases appear for chains of wavelets at multiple scales and overlapping offsets, as captured through parent--children relationships in wavelet trees \cite{BCDH10}. More general structured-sparsity models (also known as model-based sparsity and union-of-subspaces models) are often used in the literature to express these kinds of sparsity structures by allowing some supports and disallowing others~\cite{BCDH10}. The support of structured-sparse signals in many instances can be expressed in terms of unions of groups of indices $\{\Omega_1,\ldots,\Omega_M\},$ where each group $\Omega_m \subset \{1,\ldots,p\}$ contains indices of $k$ coefficients that become active simultaneously in structured-sparse signals, and the groups $\Omega_m$'s may or may not be disjoint. Therefore, despite the non-block sparse nature of $\beta$ in this setting, one can reorganize the columns of the dictionary $\Phi$ and the entries of $\beta$ as $\Phi' = [\Phi_{\Omega_1}~\Phi_{\Omega_2}~\ldots~\Phi_{\Omega_M}]$ and $\beta' = [\beta_{\Omega_1}~\beta_{\Omega_2}~\ldots~\beta_{\Omega_M}]$, respectively, so that $y = \Phi\beta = \Phi'\beta'$. (In the case of overlapping $\Omega_m$'s, this does require minor corrections to $\beta'$ to ensure each nonzero entry in $\beta$ appears only once in $\beta'$; see, e.g., \cite{RRN12,Bach10}.) Doing so converts the structured-sparse problem involving $\beta$ into a block-sparse problem involving $\beta'$ and results of Theorems~\ref{thm:rand_cond}--\ref{thm:grouplasso} can still be utilized in this case as long as the number of sets $M$ is not prohibitively large.
\end{rem}

\subsection{Numerical Experiments}
One of the most important implications of this section is that, similar to the case of recovery of block-sparse signals, the number of maximum allowable active regression variables in regression of block-sparse vectors is fundamentally independent of the column-wise coherence of the design matrix; instead, it appears to be affected by the spectral norm of the design matrix. However, such a claim needs to be carefully investigated since our results are only concerned with sufficient conditions on design matrices. To this end, we construct numerical experiments that help us evaluate the regression performance of the group lasso for a range of design matrices with varying spectral norms, coherences, inter-block coherences, intra-block coherences, and quadratic-mean block coherences. In order to generate these design matrices, we reuse the experimental setup described in Section~\ref{ssec:recovery_num_exp} (corresponding to $n = 858$, $m = 10$, and $r = 500$).

\begin{figure}
\begin{center}
\includegraphics[width=3in]{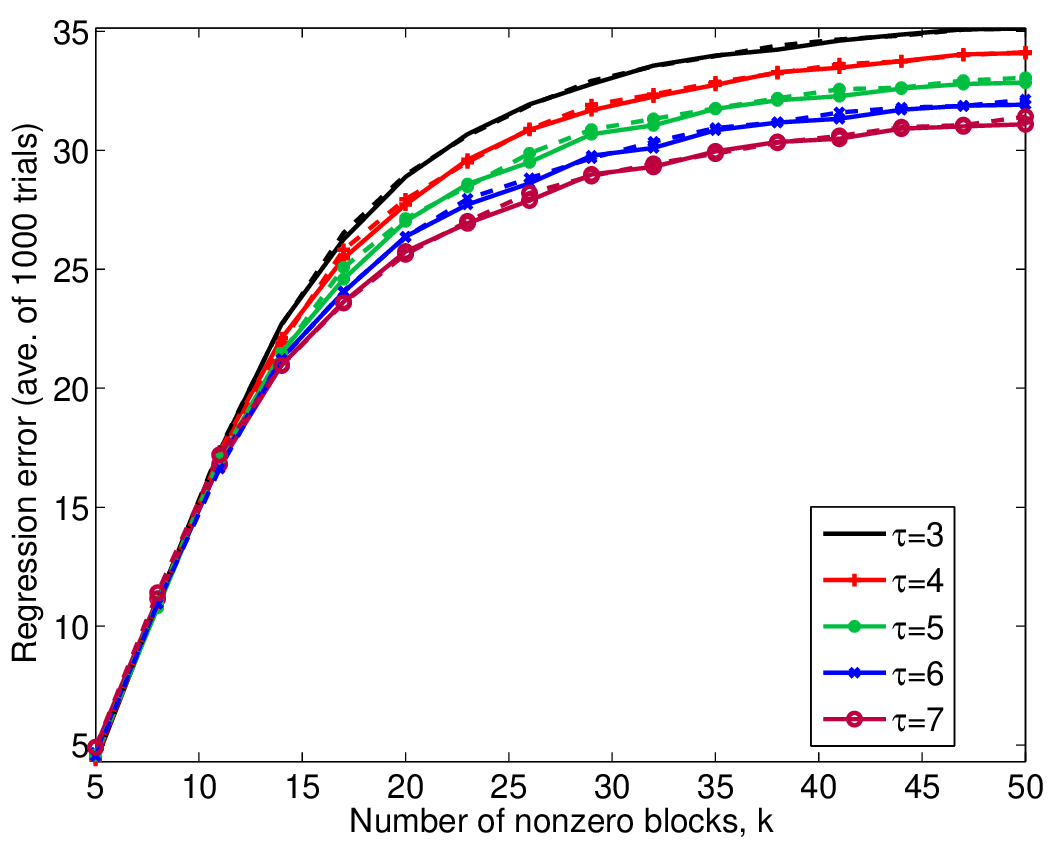}
\end{center}
\caption{Performances of the group lasso for design matrices $\Phi$ with varying spectral norms and extremal coherence values (cf. Table~\ref{table:extreme_mu}) in regression of block-sparse vectors as a function of the number of nonzero regression blocks $k$; $\tau$ denotes the value of the spectral norm multiplier used. Solid lines correspond to matrices with minimum coherence, while dashed lines correspond to matrices with maximum coherence.}
\label{fig:reg_extreme_mu}
\end{figure}

For the sake of brevity, we focus only on the performance of the group lasso \eqref{eq:grouplasso} for regression of block-sparse vectors.\footnote{We used the {\tt SpaRSA} Matlab package \cite{SPARSA} with {\tt debias} option turned on in all simulations in this section.} This performance is evaluated for different design matrices using Monte Carlo trials, corresponding to generation of 1000 block-sparse $\beta$ with $k$ nonzero blocks. Each vector of regression coefficients has block support selected uniformly at random according to M1 in Section~\ref{ssec:stat_model} and nonzero entries drawn independently from the Gaussian distribution. We then obtain the observations $y = \Phi \beta + z$ using the design matrix (dictionary) $\Phi$ under study for each one of the block-sparse $\beta$, where the variance $\sigma^2$ of the modeling error $z$ is selected such that $\|\beta\|_2^2/n\sigma^2 \approx 0.84$. Finally, we carry out linear regression using the group lasso by setting $\lambda \approx 1.4592$ and we then record the regression error $\|\Phi\beta-\Phi\betahat\|_2^2$.

Figure~\ref{fig:reg_extreme_mu} shows the regression performance of the group lasso for design matrices $\Phi$ with increasing spectral norms ($\tau = \{3,\ldots,7\}$), where we once again choose matrices with the largest and smallest coherence values for each $\tau$ (among the 2000 available options). The spectral norms, coherences, inter-block coherences, intra-block coherences, and quadratic-mean block coherences for these chosen design matrices are still given by Table~\ref{table:extreme_mu} in Section~\ref{ssec:recovery_num_exp}. Similar to the case of block-sparse recovery, we observe that significant changes in the values of the coherences do not significantly affect the regression performance. This behavior is clearly in agreement with our expectation from Theorem~\ref{thm:grouplasso} that the role of column-wise coherence in regression is decoupled from the scaling of the number of nonzero blocks $k$ (equivalently, number of nonzero regression coefficients $km$). In addition, we observe that as the spectral norm of the matrix decreases, the range of values of $k$ for which the regression error exhibits a linear trend also increases. This is again consistent with the statement of Theorem~\ref{thm:grouplasso}.

\section{Conclusion}
\label{sec:conc}
In this paper, we have provided conditions under which most block subdictionaries of a dictionary are well conditioned. In contrast to prior works, these conditions are explicitly computable in polynomial time, they can lead to near-optimal scaling of the dimensions of the well-conditioned subdictionaries for dictionaries that are approximately tight frames, and they suggest that the spectral norm plays a far important role than the column-wise coherence of the dictionary in determining the order-wise dimensions of the well-conditioned subdictionaries. In addition, we have utilized this result to investigate the average-case performances of convex optimization-based methods for block-sparse recovery and regression of block-sparse vectors. Our average-case performance results significantly improve upon the existing results for both block-sparse recovery and regression of block-sparse vectors. Finally, numerical experiments conducted in the context of block-sparse recovery and regression problems validate the insight offered by our results in relation to the effects of spectral norm and coherence of the dictionary on inference problems in block-sparse settings.

\section{Acknowledgements}
We gratefully acknowledge many helpful comments provided by Dr.~Lorne Applebaum in relation to a preliminary draft of this paper. We are also thankful to the anonymous reviewers for their many valuable remarks that helped improve the quality of this paper. Finally, we would like to thank Dr.~Jarvis Haupt for pointing out an error in the proof of the main result of our published manuscript, which led to this revision of the published version.
\appendices

\renewcommand\thelem{\thesection.\arabic{lem}}
\setcounter{lem}{0}

\section{Proof of Lemma~\ref{lem:cond_block}}
\label{sec:mainproof} The proof of this lemma relies on many lemmas and tools, some of which are generalizations of the corresponding results in~\cite{TroppCRAS,TroppRPP,Rudelson99,RudelsonVershynin} to the block setting of this paper. To begin, we denote the matrix $G$ in block-partitioned fashion:
$$G = [G_1~G_2~\ldots~G_r] = \left[\begin{array}{cccc}
G_{1,1}&G_{1,2}&\ldots&G_{1,r} \\
G_{2,1}&G_{2,2}&\ldots&G_{2,r} \\
\vdots&\vdots&\ddots&\vdots \\
G_{r,1}&G_{r,2}&\ldots&G_{r,r}
\end{array} \right],$$
where $G_{i,j} = \Phi_i^*\Phi_j-1_{\{i=j\}}\ident$ for $1 \le i,j \le r$. We then split $G = H + D$, where $D$ contains the diagonal blocks $G_{i,i}$, and $H$ contains only the non-diagonal blocks. We next define the following ``norms'' for block matrices:
\begin{itemize}
\item When we block only columns of a matrix $M$, we define $\|M\|_{B,1} := \max_{1\le i \le r} \|M_i\|_2$, and
\item When we block both columns and rows of a matrix $M$, we define $\|M\|_{B,2} := \max_{1 \le i,j \le r} \|M_{i,j}\|_2$.
\end{itemize}
Finally, we make use of some standard inequalities in the following, including:
\begin{itemize}
\item Cauchy-Schwarz Inequality: $|\E (XY)|^2 \le \E(X^2) \E(Y^2).$
\item H\"{o}lder's Inequality: $\|fg\|_1 \le \|f\|_s\|g\|_q$, $1\le s,q \le \infty$ and $1/s+1/q = 1$.
\item Jensen's Inequality for a convex function $f$: $f(\E X) \le \E f(X).$
\item Scalar Khintchine Inequality: Let $\{a_i\}$ be a finite sequence of complex numbers and $\{\epsilon_i\}$ be a Rademacher (uniformly random $\pm 1$ binary, i.i.d.) sequence. For each $q \ge 0$, we have $$\Eq \left|\sum_i \epsilon_i a_i\right| \le C_q \left(\sum_i |a_i|^2\right)^{1/2},$$ where  $C_q \le 2^{1/4} \sqrt{q/e}$.
\item Noncommutative Khintchine Inequality~\cite{TroppRPP}: Let $\{M_i\}$ be a finite sequence of matrices of the same dimensions and $\{\epsilon_i\}$ be a Rademacher sequence. For each $q \ge 2$, $$\Eq \left\|\sum_i \epsilon_i M_i \right\|_{S_q} \le W_q \max \left\{\left\|\left(\sum_i M_i M_i^*\right)^{1/2}\right\|_{S_q},\left\|\left(\sum_i M_i^* M_i\right)^{1/2}\right\|_{S_q}\right\},$$ where $\|M\|_{S_q} := \|\sigma(M)\|_q$ denotes the Schatten $q$-norm for a matrix $M$ (equal to the $\ell_q$-norm of the vector $\sigma(M)$, which contains singular values of the matrix $M$) and $W_q \le 2^{-1/4} \sqrt{\pi q/e}$.
\end{itemize}

We need the following five lemmas in our proof of Lemma~\ref{lem:cond_block}. The first two lemmas here are used to prove the later ones.
\begin{lem}
Let $X = [X_1~X_2~\ldots~X_r]$ be a block matrix and $D_X$ be its block diagonalization, i.e., a block-diagonal matrix $D_X = \mathrm{diag}(X_1,X_2, \dots, X_r)$ containing the matrices $\{X_i\}$ in its diagonal, with all other elements being equal to zero. Then, we have $$\|D_X\|_2 \le \|X\|_{B,1}.$$
\label{lem:blockdiag}
\end{lem}
\begin{IEEEproof}
For a vector $a$ of appropriate length, we evaluate the ratio $\frac{\|D_X a\|_2^2}{\|a\|_2^2}$. We partition $a = [a_1^*~a_2^*~\ldots~a_r^*]^*$ into its pieces $a_i$ matching the number of columns of the blocks $X_i$, $1 \le i \le r$. Then, we have
\begin{align*}
\frac{\|D_X a\|_2^2}{\|a\|_2^2} &= \frac{\sum_{i=1}^r \|X_i a_i\|_2^2}{\sum_{i=1}^r \|a_i\|_2^2} \le \frac{\sum_{i=1}^r \|X_i\|_2^2 \|a_i\|_2^2}{\sum_{i=1}^r \|a_i\|_2^2} \\
&\le \frac{\max_{1 \le i \le r} \|X_i\|_2^2 \sum_{i=1}^r \|a_i\|_2^2}{\sum_{i=1}^r \|a_i\|_2^2}  = \max_{1 \le i \le r} \|X_i\|_2^2.
\end{align*}
Thus, the spectral norm obeys
$$\|D_X\|_2 = \max_{a} \frac{\|D_X a\|_2}{\|a\|_2} \le \max_{1 \le i \le r} \|X_i\|_2 = \|X\|_{B,1}.$$
\end{IEEEproof}
The next lemma is a generalization of the lemma in \cite[Sec.~2]{Rudelson99} to our block setting.
\begin{lem}
Let $X = [X_1~X_2~\ldots~X_r]$ be a block matrix where each block $X_i$ has $m$ columns with $mr = p$ and let $\{\epsilon_i\}$ be a Rademacher sequence. For any $q \ge 2 \log p$, we have $$\Eq\left\|\sum_{i=1}^r \epsilon_i X_i X_i^*\right\|_2 \le 1.5\sqrt{q} \|X\|_{B,1}\|X\|_2.$$
\label{lem:rudelson}
\end{lem}
\begin{IEEEproof}
We start by bounding the spectral norm by the Schatten $q$-norm: $$E:= \Eq\left\|\sum_{i=1}^r \epsilon_i X_i X_i^*\right\|_2 \le \Eq\left\|\sum_{i=1}^r \epsilon_i X_i X_i^*\right\|_{S_q}.$$
Now we use the noncommutative Khintchine inequality (noting that the two terms in the inequality's $\max$ are equal in this case) to get $$E \le W_q \left\|\left(\sum_{i=1}^r X_i X_i^*X_iX_i^*\right)^{1/2}\right\|_{S_q}.$$ We can bound the Schatten $q$-norm by the spectral norm by paying a multiplicative penalty of $p^{1/q}$, where $p$ is the maximum rank of the matrix sum. By the hypothesis $q \ge 2 \log p$, this penalty does not exceed $\sqrt{e}$, resulting in
\begin{align*}
E &\le W_q\sqrt{e} \left\|\left(\sum_{i=1}^r X_i X_i^*X_iX_i^*\right)^{1/2}\right\|_2\\
&\le W_q\sqrt{e} \left\|\sum_{i=1}^r X_i X_i^*X_iX_i^*\right\|_2^{1/2}.
\end{align*}
Finally, we note that the sum term is a quadratic form that can be expressed in terms of $X$ and its block diagonalization, as follows:
\begin{align*}
E &\le W_q\sqrt{e} \| XD_X^*D_XX^*\|_2^{1/2} \le W_q\sqrt{e} \|D_XX^*\|_2 \\
&\le W_q\sqrt{e} \|D_X\|_2\|X\|_2 \\
&\le W_q\sqrt{e} \|X\|_{B,1}\|X\|_2,
\end{align*}
where the last step used Lemma~\ref{lem:blockdiag}. Now replace $W_q \le 2^{-1/4} \sqrt{\pi q/e}$ to complete the proof.
\end{IEEEproof}

The next lemma is a generalization of \cite[Proposition 2.1]{TroppCRAS} to our block setting.
\begin{lem}
Let $H$ be a Hermitian matrix with zero blocks on the diagonal. Then $\Eq \|RHR\|_2 \le 2 \Eq \|RHR'\|_2$, where $R' := \Sigma' \otimes \ident_m$ with $\Sigma'$ denoting an independent realization of the random matrix $\Sigma$.
\label{lem:decouple}
\end{lem}

\begin{IEEEproof}
We establish the result for $q=1$ for simplicity and without loss of generality. Denote by $\Hhat_{i,j}$ the masking of the matrix $H$ that preserves only the subblock $H_{i,j}$ and makes other entries of $H$ zero. Then, we have $$\E\|RHR\|_2 = \E\left\|\sum_{1 \le i < j \le r} \zeta_i\zeta_j (\Hhat_{i,j}+\Hhat_{j,i})\right\|_2.$$ Let $\eta_i$ be i.i.d.\ Bernoulli random variables with parameter $1/2$. We then use Jensen's inequality on this new (multivariate) random variable $\eta = \{\eta_i\}_{1 \le i \le r}$. Specifically, we define $M_{i,j}(\eta)=\eta_i(1-\eta_j)+\eta_j(1-\eta_i)$, and note that $\E_\eta M_{i,j}(\eta) = 1/2$ for all $i,j$, where $\E_\eta$ denotes expectation over the random variable $\eta$ (in contrast to the notation $\E_q$, where $q$ is a constant). We also define the function $$f(M_{i,j}(\eta)) =  \E_\zeta \left\|\sum_{1 \le i < j \le r} 2\zeta_i\zeta_jM_{i,j}(\eta)(\Hhat_{i,j}+\Hhat_{j,i})\right\|_2.$$ Then, by applying Jensen's inequality to $f(\E_\eta M_{i,j}(\eta))=\E\|RHR\|_2$, we obtain $$\E\|RHR\|_2 \le 2 \E_\eta \E_\zeta \left\|\sum_{1 \le i < j \le r} [\eta_i(1-\eta_j)+\eta_j(1-\eta_i)]\zeta_i\zeta_j (\Hhat_{i,j}+\Hhat_{j,i})\right\|_2.$$ There is a 0-1 vector $\eta^*$ for which the expression exceeds its expectation over $\eta$. Letting $T = \{i : \eta_i^* = 1\}$, we get
\begin{align}
\nonumber
\E\|RHR\|_2 &\le 2\E_\zeta \left\|\sum_{i \in T, j \in T^C} \zeta_i\zeta_j (\Hhat_{i,j}+\Hhat_{j,i})\right\|_2 = 2\E_\zeta \left\|\sum_{i \in T, j \in T^C} \zeta_i\zeta_j \Hhat_{i,j} + \sum_{i \in T, j \in T^C} \zeta_i\zeta_j \Hhat_{j,i}\right\|_2\\
 &= 2 \E_\zeta \left\|\sum_{i \in T, j \in T^C} \zeta_i\zeta_j \Hhat_{i,j}\right\|_2 = 2 \E_\zeta \left\|\sum_{i \in T, j \in T^C} \zeta_i\zeta_j' \Hhat_{i,j}\right\|_2,\label{eq:dilation}
\end{align}
where $\{\zeta_j'\}$ is an independent realization of the sequence $\{\zeta_i\}$. The equality in \eqref{eq:dilation} is a combination of the following facts: ($i$) $\Hhat_{i,j} = \Hhat_{j,i}^*$ and, therefore, defining $A := \sum_{i \in T, j \in T^C} \zeta_i\zeta_j \Hhat_{i,j}$ and $B := \sum_{i \in T, j \in T^C} \zeta_i\zeta_j \Hhat_{j,i}$, we have $B = A^*$; ($ii$) we can reorder the columns and rows of $A$ and $B$ to have those corresponding to the set $T$ first, followed by those corresponding to the set $T^C$ later, giving us matrices of the form $$\check{A} = \begin{bmatrix}0 & \tilde{A}\\ 0 & 0\end{bmatrix}~\textrm{and}~\check{B} = \begin{bmatrix}0 & 0\\ \tilde{B} & 0\end{bmatrix},$$ respectively, where $\tilde{A} = \tilde{B}^*$; ($iii$) permuting the columns and rows of a matrix does not affect its spectral norm; ($iv$) the Hermitian dilation map of a matrix $M$, $$\mathcal{D}(M) : M \mapsto \begin{bmatrix}0 & M\\ M^* & 0\end{bmatrix},$$ preserves the spectral norm of $M$: $\|M\|_2 = \|\mathcal{D}(M)\|_2$; and ($v$) removal of all-zero rows and columns from a matrix preserves its spectral norm. Combining these five facts together, we have $$\|A+B\|_2 = \|\check{A}+\check{B}\|_2 = \left\|\begin{bmatrix}0 & \tilde{A}\\ \tilde{B}& 0\end{bmatrix}\right\|_2 = \|\mathcal{D}(\tilde{A})\|_2 = \|\tilde{A}\|_2 = \|\check{A}\|_2 = \|A\|_2.$$ Finally, since the norm of a submatrix does not exceed the norm of the matrix, we re-introduce the missing blocks to complete the argument: $$\E\|RHR\|_2 \le  2 \E_\zeta \left\|\sum_{1 \le i,j \le r, i \ne j} \zeta_i\zeta_j' \Hhat_{i,j}\right\|_2 = 2 \E_\zeta \|RHR'\|_2.$$
\end{IEEEproof}
The next lemma is adapted to our problem setup of block matrices from~\cite[Theorem~3.1]{RudelsonVershynin}, \cite[Proposition~12]{TroppRPP}.
\begin{lem}
Let $M = [M_1~M_2~\ldots~M_r]$ be a matrix with $r$ column blocks, and suppose $q \ge 4 \log p \ge 2$. Then $$\Eq \|MR\|_2 \le 3\sqrt{\frac{q}{2}}\Eq\|MR\|_{B,1} + \sqrt{\delta} \|M\|_2.$$
\label{lem:RV}
\end{lem}
\begin{IEEEproof}
We denote $E := \Eq \|MR\|_2$ and note that
\begin{align*}
E^2 &= \E_{q/2}\|MRM^*\|_2 = \E_{q/2} \left\|\sum_{1\le i \le r} \zeta_i M_i M_i^*\right\|_2 \\
&\le \E_{q/2} \left\|\sum_{1\le i \le r} (\zeta_i-\delta) M_i M_i^*\right\|_2 + \delta\left\|\sum_{1\le i \le r} M_i M_i^*\right\|_2.
\end{align*}
Next, we replace $\delta$ by $\delta = \E \zeta_i'$, with $\{\zeta_i'\}$ denoting an independent copy of the sequence $\{\zeta_i\}$. We then take the expectation out of the norm by applying Jensen's inequality to get $$E^2 \le \E_{q/2} \left\|\sum_{1\le i \le r} (\zeta_i-\zeta_i') M_i M_i^*\right\|_2 + \delta\|MM^*\|_2.$$ We now symmetrize the distribution by introducing a Rademacher sequence $\{\epsilon_i\}$, noticing that the expectation does not change due to the symmetry of the random variables $\zeta_i-\zeta_i'$: $$E^2 \le \E_{q/2} \left\|\sum_{1\le i \le r} \epsilon_i(\zeta_i-\zeta_i') M_i M_i^*\right\|_2 + \delta\|M\|_2^2.$$ We apply the triangle inequality to separate $\zeta_i$ and $\zeta_i'$, and by noticing that they have the same distribution, we obtain $$E^2 \le 2\E_{q/2} \left\|\sum_{1\le i \le r} \epsilon_i\zeta_i M_i M_i^*\right\|_2 + \delta\|M\|_2^2.$$ Writing $\Omega = \{i : \zeta_i = 1\}$, we see that $$E^2 \le 2\E_{q/2,\zeta}\left(\E_{q/2,\epsilon} \left\|\sum_{i \in \Omega} \epsilon_i M_i M_i^*\right\|_2\right) + \delta\|M\|_2^2,$$ where we have split the expectation on the random variables $\{\zeta_i\}$ and $\{\epsilon_i\}$. Now we use Lemma~\ref{lem:rudelson} on the term in parentheses to get $$E^2 \le 3 \sqrt{\frac{q}{2}} \E_{q/2} (\|MR\|_{B,1}\|MR\|_2)  + \delta\|M\|_2^2.$$ Using the Cauchy-Schwarz inequality, we get $$E^2 \le 3 \sqrt{\frac{q}{2}} \Eq \|MR\|_{B,1} \Eq \|MR\|_2  + \delta\|M\|_2^2.$$ This inequality takes the form $E^2 \le bE+c$. We bound $E$ by the largest solution of this quadratic form: $$E \le \frac{b+\sqrt{b^2+4c}}{2} \le b+\sqrt{c},$$ thereby proving the lemma.
\end{IEEEproof}
The last lemma that we need for our proof is a generalization of \cite[Proposition 13]{TroppRPP} to our block setting.
\begin{lem}
Let
$$M = [M_1~M_2~\ldots~M_r] = \left[\begin{array}{cccc}
M_{1,1}&M_{1,2}&\ldots&M_{1,r} \\
M_{2,1}&M_{2,2}&\ldots&M_{2,r} \\
\vdots&\vdots&\ddots&\vdots \\
M_{r,1}&M_{r,2}&\ldots&M_{r,r}
\end{array} \right]$$
be a block matrix, where each block $M_{i,j}$ has size $m \times m$. Assume $q \ge 2 \log r$.
Then, we have
$$\Eq\|RM\|_{B,1} \le 2^{1.5}\sqrt{q} \|M\|_{B,2} +\sqrt{\delta \max_{1 \le j \le r} \sum_{i=1}^r \|M_{i,j}\|_2^2}.$$
\label{lem:Bnorm}
\end{lem}
\begin{IEEEproof}
We begin by seeing that
$$E^2 := (\Eq \|RM\|_{B,1})^2 = \left[\Eq\left( \max_{1 \le j \le r} \|RM_j\|_2\right)\right]^2 \le \E_{q/2} \left(\max_{1 \le j \le r} \sum_{i=1}^r \zeta_i \|M_{i,j}\|_2^2\right),$$ where the inequality can be derived using the definition of the spectral norm. In the sequel, we abbreviate $t = q/2$, $y_{i,j} = \|M_{i,j}\|_2^2$, and $\E_t \left(\max_j \sum_i \zeta_i y_{i,j} \right) = \widetilde{E}^2$. We now bound $\widetilde{E}^2$ by continuing to use the same technique as in the proof of Lemma~\ref{lem:RV}: we split a term for the mean value of the sequence $\{\zeta_i\}$, then replace the term by $\{\E \zeta_i'\}$, with $\{\zeta_i'\}$ an independent copy of the sequence $\{\zeta_i\}$, then apply Jensen's inequality and exploit symmetrization by introducing a Rademacher sequence $\{\epsilon_i\}$, and finally finish by merging the two terms due to their identical distributions:
\begin{align*}
\widetilde{E}^2 &\le \E_{t} \left(\max_{1 \le j \le r} \sum_{i=1}^r (\zeta_i-\delta) y_{i,j}\right)+\delta \max_{1 \le j \le r} \sum_{i=1}^r y_{i,j} \\
&\le \E_{t} \left(\max_{1 \le j \le r} \sum_{i=1}^r (\zeta_i-\zeta_i') y_{i,j}\right)+\delta \max_{1 \le j \le r} \sum_{i=1}^r y_{i,j} \\
&= \E_{t} \left(\max_{1 \le j \le r} \sum_{i=1}^r \epsilon_i(\zeta_i-\zeta_i') y_{i,j}\right)+\delta \max_{1 \le j \le r} \sum_{i=1}^r y_{i,j} \\
&\le 2\E_{t} \left(\max_{1 \le j \le r} \sum_{i=1}^r \epsilon_i\zeta_i y_{i,j}\right)+\delta \max_{1 \le j \le r} \sum_{i=1}^r y_{i,j}.
\end{align*}
Now we bound the maximum by the $\ell_t$ norm and separate the expectations on the two sequences:
$$\widetilde{E}^2 \le 2 \left(\E_{\zeta}\sum_{j=1}^r \left(\E_{t,\epsilon} \sum_{i=1}^r \epsilon_i \zeta_i y_{i,j}\right)^{t}\right)^{1/t}+\delta \max_{1 \le j \le r} \sum_{i=1}^r y_{i,j}.$$
For the inner term, we can use the scalar Khintchine inequality to obtain
$$\widetilde{E}^2 \le 2C_t \left(\E_{\zeta} \sum_{j=1}^r \left(\sum_{i=1}^r \zeta_i y_{i,j}^2\right)^{t/2}\right)^{1/t}+\delta \max_{1 \le j \le r} \sum_{i=1}^r y_{i,j}.$$
We continue by bounding the outer sum by the maximum term times the number of terms:
$$\widetilde{E}^2 \le 2C_t r^{1/t} \left(\E_{\zeta} \max_{1 \le j \le r} \left(\sum_{i=1}^r \zeta_i y_{i,j}^2\right)^{t/2}\right)^{1/t}+\delta \max_{1 \le j \le r} \sum_{i=1}^r y_{i,j}.$$
Since $t \ge \log r$, it holds that $r^{1/t} \le e$, which implies that
$2C_tr^{1/t} \le 4\sqrt{t}$.  We now use H\"older's inequality inside the sum
term $\zeta_i y_{i,j}^2 = y_{i,j}\cdot \zeta_i y_{i,j}$ with $s=\infty$,
$q=1$:
\begin{align*}
\widetilde{E}^2 &\le 4\sqrt{t}\left(\max_{1\le i,j \le r} y_{i,j}\right)^{1/2} \left(\E_{\zeta} \max_{1 \le j \le r} \left(\sum_{i=1}^r \zeta_i y_{i,j}\right)^{t/2}\right)^{1/t}+\delta \max_{1 \le j \le r} \sum_{i=1}^r y_{i,j} \\
&\le 4\sqrt{t}\left(\max_{1\le i,j \le r} y_{i,j}\right)^{1/2} \left(\E_{\zeta} \max_{1 \le j \le r} \left(\sum_{i=1}^r \zeta_i y_{i,j}\right)^{t}\right)^{1/2t}+\delta \max_{1 \le j \le r} \sum_{i=1}^r y_{i,j},
\end{align*}
where the last inequality follows from the concavity of the square root. It then follows that
\begin{align*}
\widetilde{E}^2 &\le 4\sqrt{t} \max_{1\le i,j \le r} \sqrt{y_{i,j}} \sqrt{\E_{t,\zeta} \left(\max_{1 \le j \le r} \sum_{i=1}^r \zeta_i y_{i,j}\right)}+\delta \max_{1 \le j \le r} \sum_{i=1}^r y_{i,j}\\
&= 4\sqrt{t} \max_{1\le i,j \le r} \sqrt{y_{i,j}}\widetilde{E} +\delta \max_{1 \le j \le r} \sum_{i=1}^r y_{i,j}.
\end{align*}
Given that $\widetilde{E}$ has also appeared on the right hand side, we follow the same argument that ends the proof of Lemma~\ref{lem:RV} to obtain
\begin{align*}
\widetilde{E} &\le 4\sqrt{t} \max_{1\le i,j \le r} \sqrt{y_{i,j}} + \sqrt{\delta \max_{1 \le j \le r} \sum_{i=1}^r y_{i,j}}.
\end{align*}
Now we recall that $E \leq \widetilde{E}$, $t = q/2$ and $y_{i,j} = \|M_{i,j}\|_2^2$ to get
\begin{align*}
E &\le 2^{1.5}\sqrt{q}\max_{1\le i,j \le r} \|M_{i,j}\|_2 + \sqrt{\delta \max_{1 \le j \le r} \sum_{i=1}^r \|M_{i,j}\|_2^2}\\
&= 2^{1.5}\sqrt{q}\|M\|_{B,2} + \sqrt{\delta \max_{1 \le j \le r} \sum_{i=1}^r \|M_{i,j}\|_2^2}.
\end{align*}
This completes the proof of the lemma.
\end{IEEEproof}

We now have all the required results to prove Lemma~\ref{lem:cond_block}. Split $G$ into its diagonal blocks $D$ (containing $\Phi_i^*\Phi_i-\ident$, $1\le i \le r$) and off-diagonal blocks $H$ (containing $\Phi_i^*\Phi_j$, $1 \le i\ne j \le r$) and apply Lemma~\ref{lem:decouple}: $$\Eq \|RGR\|_2 \le 2 \Eq \|RHR'\|_2 + \Eq \|RDR\|_2.$$ To estimate the first term, we apply Lemma~\ref{lem:RV} twice; once for $R$, and once for $R'$:
\begin{align*}
\Eq \|RHR'\|_2 &\le 3\sqrt{\frac{q}{2}} \Eq \|RHR'\|_{B,1} + \sqrt{\delta} \Eq \|HR'\|_2\\
&\le 3\sqrt{\frac{q}{2}} \Eq \|RHR'\|_{B,1} + 3\sqrt{\frac{\delta q}{2}} \Eq \|HR'\|_{B,1} + \delta \|H\|_2.
\end{align*}
By applying Lemma~\ref{lem:Bnorm} on the first term, we obtain
\begin{align*}
\Eq \|RHR'\|_2 \le 3\sqrt{\frac{q}{2}}\left(\sqrt{8q}\Eq \|HR'\|_{B,2}+\sqrt{\delta \max_{1 \le j \le r} \sum_{i=1}^r \|(HR')_{i,j}\|_2^2}\right) + 3\sqrt{\frac{\delta q}{2}} \Eq \|HR'\|_{B,1} + \delta \|H\|_2.
\end{align*}
Since $R$ and $R'$ have the same distribution, we get $$ \Eq \|RGR\|_2 \le 12q\Eq \|HR\|_{B,2}+6\sqrt{\frac{\delta q}{2} \max_{1 \le j \le r} \sum_{i=1}^r \|(HR)_{i,j}\|_2^2} + 6\sqrt{\frac{\delta q}{2}} \Eq \|HR\|_{B,1} + 2\delta \|H\|_2 + \Eq \|RDR\|_2.$$ Next, in order to bound $\|HR\|_{B,1}$, we use the notation $\Phi_{\{i\}^C} = [\Phi_1~\ldots~\Phi_{i-1}~\Phi_{i+1}~\ldots~\Phi_r]$; we then have
\begin{align*}
\|HR\|_{B,1} &\le \|H\|_{B,1} = \max_{1 \le i \le r} \|\Phi_i^*\Phi_{\{i\}^C}\|_2 \le \max_{1 \le i \le r} \|\Phi_i^*\Phi\|_2 \\
&\le \max_{1 \le i \le r} \|\Phi_i\|_2\|\Phi\|_2 = \sqrt{1+\mu_I}\|\Phi\|_2.
\end{align*}
Now we use the facts $\|HR\|_{B,2} \le \mu_B$, $\sqrt{\max_{1 \le j \le r} \sum_{i=1}^r \|(HR)_{i,j}\|_2^2} \leq \sqrt{r-1}\,\overline{\mu}_B$, $\|H\|_2 \le \|G\|_2+\|D\|_2 = \|\Phi\|_2^2+\|D\|_2$ and, using Lemma~\ref{lem:blockdiag}, $$\Eq \|RDR\|_2 \le \|D\|_2 \leq \max_{1 \le i \le r} \|\Phi_i^*\Phi_i -\ident\|_2 \le \mu_I$$ to complete the proof of the lemma:
\begin{align*}
\Eq \|RGR\|_2 &\le 12q\mu_B +6\sqrt{\frac{\delta q(r-1)}{2}}\overline{\mu}_B + 6\sqrt{\frac{\delta q(1+\mu_I)}{2}}\|\Phi\|_2 + 2\delta(\|\Phi\|_2^2+\mu_I)+\mu_I \\
&\le 48\mu_B\log p +6\overline{\mu}_B\sqrt{2(r-1)\delta \log p} + 6\sqrt{2\delta\log{p}(1+\mu_I)}\|\Phi\|_2 + 2\delta\|\Phi\|_2^2+3\mu_I.
\end{align*}

\setcounter{lem}{0}

\section{Proof of Theorem~\ref{thm:l12}}
\label{sec:l12proof} In this appendix, we will prove that the minimization (\ref{eq:l21min}) successfully recovers a $k$-block sparse $\beta$ from $y = \Phi \beta$ with high probability. Mathematically, this is equivalent to showing that $\|\beta\|_{2,1} < \|\beta'\|_{2,1}$ for all $\beta' \ne \beta$ such that $y = \Phi \beta'$. In the following, we will argue that this is true as long as there exists a vector $h \in \real^{km}$ such that ($i$) $\Phi_\cl{S}^* h = \signbar(\beta_\cl{S})$, where $\cl{S}$ denotes the block support of $\beta$, $\Phi_\cl{S}^*$ denotes the adjoint of $\Phi_\cl{S}$ (rather than a column submatrix of $\Phi^*$), and $\signbar(\beta_\cl{S})$ denotes the block-wise extension of $\signbar(\cdot)$ to the blocks in $\cl{S}$, and ($ii$) $\|\Phi_j^*h\|_2 < 1$ for all $j \notin \cl{S}$. Note that these two conditions on the vector $h$ imply that ($iii$) $\|\Phi_j^*h\|_2 \le 1$ for all $1\le j \le r$.

To prove the sufficiency of conditions ($i$) and ($ii$) above, we follow the same ideas as in~\cite{Fuchs04,TroppIT05,ER10}. To begin, we need the following lemma; its proof is a simple exercise using H\"{o}lder's Inequality.
\begin{lem}
Consider two block vectors $a$ and $b$ such that the blocks of $a$ have nonidentical $\ell_2$-norms and the blocks of $b$ are nonzero. Then $\langle a,b \rangle < \|a\|_{2,\infty} \|b\|_{2,1}$, where $\|a\|_{2,\infty} := \max_{j=1,\ldots,r} \|a_j\|_2$. \label{lem:holder2}
\end{lem}
\begin{IEEEproof}
We can write
\begin{align}
\langle a,b \rangle = \sum_{j=1}^r\sum_{l=1}^m a_{jl}b_{jl} = \sum_{j=1}^r\langle a_j,b_j\rangle \le \sum_{j=1}^r\|a_j\|_2\|b_j\|_2 = \langle\overline{a},\overline{b}\rangle,
\label{eq:blockcs}
\end{align}
where the vectors $\overline{a} = [\|a_1\|_2~\|a_2\|_2~\ldots~\|a_r\|_2]$ and $\overline{b} = [\|b_1\|_2~\|b_2\|_2~\ldots~\|b_r\|_2]$ are defined as ones that contain the $\ell_2$-norms of the blocks of $a$ and $b$, respectively. Using the conditions of~\cite[Lemma 6]{TroppIT05}, we have $$\langle a,b \rangle \le \langle\overline{a},\overline{b}\rangle < \|\overline{a}\|_\infty \|\overline{b}\|_1 = \|a\|_{2,\infty} \|b\|_{2,1},$$ thereby proving the lemma.
\end{IEEEproof}

\begin{rem}
Note that if we remove the requirements on $a$ and $b$, it can be shown that $\langle a,b \rangle \le \|a\|_{2,\infty} \|b\|_{2,1}$; that is, the conditions on $a$ and $b$ remove the possibility of equality.
\end{rem}

We are ready now to formally prove Theorem~\ref{thm:l12}. We begin by writing
\begin{align*}
\|\beta\|_{2,1}&=\|\beta_\cl{S}\|_{2,1}=\sum_{j \in \cl{S}} \|\beta_j\|_2 = \sum_{j \in \cl{S}} \frac{\beta_j^*\beta_j}{\|\beta_j\|_2}
 = \sum_{j \in \cl{S}} \left\langle \frac{\beta_j}{\|\beta_j\|_2},\beta_j\right\rangle\\
& = \sum_{j \in \cl{S}} \left\langle \signbar(\beta_j),\beta_j\right\rangle = \left\langle \signbar(\beta_\cl{S}),\beta_\cl{S}\right\rangle.
\end{align*}
Next, we assume that conditions ($i$) and ($ii$) (which together imply condition ($iii$)) are true in our case and consider any $\beta' \not= \beta$ such that $y=\Phi\beta'$. Then since $\signbar(\beta_\cl{S}) = \Phi_\cl{S}^*h$, we have $$\|\beta\|_{2,1} = \left\langle \Phi_\cl{S}^*h,\beta_\cl{S}\right\rangle = \left\langle h,\Phi_\cl{S}\beta_\cl{S}\right\rangle = \langle h,y\rangle=\left\langle h,\Phi_{\cl{S}'}\beta'_{\cl{S}'}\right\rangle = \left\langle \Phi_{\cl{S}'}^*h,\beta'_{\cl{S}'}\right\rangle,$$ where $\cl{S}'$ denotes the support of a different solution $\beta'$ as described earlier. We now consider two cases. If not all norms $\|\Phi_j^*h\|_2$ are identical over $j \in \cl{S}'$, then we apply Lemma~\ref{lem:holder2} to obtain $$\|\beta\|_{2,1} < \|\Phi_{\cl{S}'}^*h\|_{2,\infty}\|\beta'_{\cl{S}'}\|_{2,1} \le \|\beta'_{\cl{S}'}\|_{2,1} =  \|\beta'\|_{2,1},$$ where the last inequality is due to condition ($iii$). If all the norms $\|\Phi_j^*h\|_2$ are identical over $j \in \cl{S}'$, note that since $\beta \ne \beta'$ and since Theorem~\ref{thm:rand_cond} guarantees that $\Phi_\cl{S}$ has linearly independent columns with high probability (noting that we will come back to Theorem~\ref{thm:rand_cond} later), then there must exist a block index $j_0 \in \cl{S}'$ such that $j_0 \notin \cl{S}$. From condition ($ii$), we know that for such a $j_0$ we have $\|\Phi_{j_0}^*h\|_2 < 1$, meaning that $\|\Phi_j^*h\|_2 < 1$ for all $j \in \cl{S}'$. We then leverage (\ref{eq:blockcs}) to obtain $$\|\beta\|_{2,1} \le \sum_{j\in \cl{S}'}\|\Phi_j^*h\|_2\|\beta'_j\|_2 = \|\Phi_{j_0}^*h\|_2\sum_{j\in \cl{S}'}\|\beta'_j\|_2 < \|\beta'_{\cl{S}'}\|_{2,1} = \|\beta'\|_{2,1}.$$ In order to complete the proof of the theorem, the only thing that remains to be shown now is that conditions ($i$) and ($ii$) hold in our case.

To simplify conditions ($i$) and ($ii$), we can define the vector $h = (\Phi_\cl{S}\pinv)^*\signbar(\beta_\cl{S})$, where $(\cdot)\pinv$ denotes the Moore--Penrose pseudoinverse of a matrix. Note that such an $h$ trivially satisfies condition ($i$). Condition ($ii$) then reduces to
\begin{align}
\|\Phi_{\cl{S}^C}^*(\Phi_\cl{S}\pinv)^*\signbar(\beta_\cl{S})\|_{2,\infty} < 1.
\label{eq:l12cond}
\end{align}
This inequality is implied by
\begin{align}
\|\Phi_{\cl{S}^C}^*(\Phi_\cl{S}\pinv)^*\signbar(\beta_\cl{S})\|_{\infty} < 1/\sqrt{m}.
\label{eq:linfcond}
\end{align}
It then remains to prove (\ref{eq:linfcond}). To this end, we condition on the event that $\Phi_\cl{S}$ has linearly independent columns and denote the scalar $$Z_{0,ij} = \phi_{i,j}^*\Phi_\cl{S}(\Phi_\cl{S}^*\Phi_\cl{S})^{-1}\signbar(\beta_\cl{S})$$ for each $i \in \cl{S}^C$ and $j = 1,\ldots,m$, where $\phi_{i,j}$ denotes the $j^{th}$ column of the $i^{th}$ block of $\Phi$. Further, we define $$Z_0 = \max_{i \notin \cl{S}, j=1,\ldots,m}|Z_{0,ij}| = \|\Phi_{\cl{S}^C}^*\Phi_\cl{S}(\Phi_\cl{S}^*\Phi_\cl{S})^{-1}\signbar(\beta_\cl{S})\|_{\infty}.$$
We now establish that $Z_0 < 1/\sqrt{m}$ with high probability. For convenience, define $W_{i,j} = (\Phi_\cl{S}^*\Phi_\cl{S})^{-1}\Phi_\cl{S}^*\phi_{i,j}$; we can then write $Z_{0,ij} = \langle W_{i,j},\signbar(\beta)\rangle = \sum_{l \in \cl{S}}\langle W_{ij,l},\signbar(\beta_l)\rangle$, where $W_{ij,l}$ denotes the block of the vector $W_{i,j}$ corresponding to the entry $l \in \cl{S}$. It then follows that
\begin{align*}
|\langle W_{ij,l},\signbar(\beta_l)\rangle| \le \|W_{ij,l}\|_2 \|\signbar(\beta_l)\|_2 = \|W_{ij,l}\|_2.
\end{align*}
Now, we use Hoeffding's inequality (similar to~\cite[Lemma 3.3]{CandesPlan}) to obtain
\begin{align*}
\prob(|Z_{0,ij}| \geq t) \le 2e^{-t^2/2\sum_{l \in \cl{S}}\|W_{ij,l}\|_2^2} = 2e^{-t^2/2\|W_{i,j}\|_2^2}.
\end{align*}
A union bound then gives us $\prob(Z_0 \ge t) \le 2pe^{-t^2/2\kappa^2}$, where $\kappa \geq \max_{i \notin \cl{S}, 1\le j \le m}\|W_{i,j}\|_2$. We can then condition on the event $\|(\Phi_\cl{S}^*\Phi_\cl{S})^{-1}\|_2 \leq 2$ and obtain
\begin{align*}
\max_{i \notin \cl{S}, 1 \le j \le m}\|W_{i,j}\|_2 &= \max_{i \notin \cl{S}, 1\le j \le m}\|(\Phi_\cl{S}^*\Phi_\cl{S})^{-1}\Phi_\cl{S}^*\phi_{i,j}\|_2 \le 2\max_{i \notin \cl{S}, 1 \le j \le m}\|\Phi_\cl{S}^*\phi_{i,j}\|_2\nonumber \\
&\le 2\max_{i \notin \cl{S}}\|\Phi_\cl{S}^*\Phi_i\|_2 = 2\|\Phi_\cl{S}^*\Phi_{\cl{S}^C}\|_{B,1}.
\end{align*}
Thus, conditioned on the bounds
\begin{align}
\|\Phi_\cl{S}^*\Phi_{\cl{S}^C}\|_{B,1} \le \gamma
\label{eq:gamma1}
\end{align}
and $\|(\Phi_\cl{S}^*\Phi_\cl{S})^{-1}\|_2 \leq 2$, and replacing $t=1/\sqrt{m}$, the probability of the inequality (\ref{eq:linfcond}) failing to hold is at most $2pe^{-1/8m\gamma^2}$.

To finalize, we define $Z = \|\Phi_\cl{S}^*\Phi_\cl{S}-\ident\|_2$ and define the event $$E = \{Z \le 1/2\} \cap \{\|\Phi_\cl{S}^*\Phi_{\cl{S}^C}\|_{B,1} \le \gamma\}.$$ Then we have that the probability $P$ of the condition (\ref{eq:l12cond}) not being satisfied is upper bounded by
\begin{align}
\label{eqn:proof_thm2_prob}
P &\le \prob(\{Z_0 \geq 1/\sqrt{m}\}|E)+\prob(E^C) \le 2 pe^{-1/8 m\gamma^2} + \prob(Z > 1/2) + \prob(\|\Phi_\cl{S}^*\Phi_{\cl{S}^C}\|_{B,1} > \gamma).
\end{align}
We set $\gamma = c_3/\sqrt{m\log p}$ with small enough $c_3$ so that the first term of the right hand side is upper bounded by $2 p^{-4\log 2}$.  Next, we appeal to Lemma~\ref{lem:Bnorm} together with the Markov inequality and a Poissonization argument (see (\ref{eq:poissonization}) and (\ref{eq:poisscp_prf}) for an example) to obtain
\begin{align}
\prob(\|\Phi_\cl{S}^*\Phi_{\cl{S}^C}\|_{B,1} > \gamma) &\le 2 \prob(\|R H\|_{B,1} > \gamma) \le 2\gamma^{-q}\mathbb{E}(\|R H\|_{B,1}^q) \nonumber\\
&\le 2\gamma^{-q}(2^{1.5}\sqrt{q}\mu_B+\sqrt{\delta(r-1)}\overline{\mu}_B)^q\nonumber \\
&\le 2\gamma^{-q}(2^{1.5}\sqrt{q}\mu_B+\sqrt{k}\overline{\mu}_B)^q,
\label{eqn:Appendix.B.end}
\end{align}
where $q = 4\log p$. We replace the values of $\gamma$ and $q$ selected above as well as the bounds on $\mu_B$ and $k$ from the theorem statement to obtain
\begin{align*}
\prob\left(\|\Phi_\cl{S}^*\Phi_{\cl{S}^C}\|_{B,1} > \frac{c_3}{\sqrt{m\log{p}}}\right) &\le 2\left(4\sqrt{2}\frac{c_2'}{c_3} + \frac{\sqrt{c_0''}}{c_3}\right)^{4\log{p}}.
\end{align*}
By picking the constants $c_0'',c_2'$ small enough so that the base of the exponential term on the right hand side is less than $1/2$, we get $\prob(\|\Phi_\cl{S}^*\Phi_{\cl{S}^C}\|_{B,1} > \frac{c_3}{\sqrt{m\log{p}}}) < 2p^{-4 \log 2}$. It therefore follows from \eqref{eqn:proof_thm2_prob} and Theorem~\ref{thm:rand_cond} (since the theorem conditions on $\mu_I$ and $\mu_B$ imply the BIC) that \eqref{eq:l12cond} holds with probability at least $1-6p^{-4\log2}$. This completes the proof of the theorem.

\setcounter{lem}{0}

\section{Proof of Theorem~\ref{thm:lasso}}
\label{app:thm:lasso}
Similar to the proof in~\cite{CandesPlan} for linear regression in the non-block setting, the proof of this theorem relies on three conditions involving the design matrix $\Phi$, the vector of regression coefficients $\beta$, and the modeling error $z$. We once again use $\cl{S}$ to denote the block support of the $k$-block sparse $\beta$, $\Phi_\cl{S}$ to denote the matrix containing columns of the blocks indexed by $\cl{S}$, i.e., an $n \times km$ submatrix of $\Phi$, and $\Phi_\cl{S}^*$ to denote the adjoint of $\Phi_\cl{S}$. Finally, we assume in this appendix that $\sigma = 1$ without loss of generality. We then have from the analysis for the lasso in \cite{CandesPlan} that the following three conditions are sufficient for the theorem statement to hold:
\begin{itemize}
\item {\em Invertibility condition}: The submatrix $\Phi_\cl{S}^*\Phi_\cl{S}$ is invertible and obeys $\|(\Phi_\cl{S}^*\Phi_\cl{S})^{-1}\|_2 \le 2$.

\item {\em Orthogonality condition}: The vector $z$ obeys $\|\Phi^*z\|_\infty \le \sqrt{2}\lambda$.

\item {\em Complementary size condition}: The following inequality holds: $$\|\Phi_{\cl{S}^C}^*\Phi_\cl{S}(\Phi_\cl{S}^*\Phi_\cl{S})^{-1}\Phi_\cl{S}^*z\|_\infty +2\lambda\|\Phi_{\cl{S}^C}^*\Phi_\cl{S}(\Phi_\cl{S}^*\Phi_\cl{S})^{-1}\sign(\beta_\cl{S})\|_\infty \le (2-\sqrt{2})\lambda.$$
\end{itemize}

In order to prove this theorem, therefore, we need only evaluate the probability with which each condition holds under the assumed statistical model, after which a simple union bound will get us the desired result. In this regard, we already have from Theorem~\ref{thm:rand_cond} that the invertibility condition fails to hold with probability at most $2 p^{-4\log 2}$. It is also easy to see that the orthogonality condition in our case is the same as that in~\cite{CandesPlan}, where it is shown that the condition fails to hold with probability at most $p^{-1}(2\pi \log p)^{-1/2}$. Therefore, we focus only on understanding the probability of failure for the complementary size condition in the following.

We begin by posing two separate statements that imply the condition, following~\cite[p. 2167]{CandesPlan}:
\begin{align}
\|\Phi_{\cl{S}^C}^*\Phi_\cl{S}(\Phi_\cl{S}^*\Phi_\cl{S})^{-1}\sign(\beta_\cl{S})\|_\infty &\le 1/4, \label{eq:csc1}\\
\|\Phi_{\cl{S}^C}\Phi_\cl{S}(\Phi_\cl{S}^*\Phi_\cl{S})^{-1}\Phi_\cl{S}^*z\|_\infty &\le (3/2-\sqrt{2})\lambda.  \label{eq:csc2_1}
\end{align}
First, we consider the inequality (\ref{eq:csc1}). Denote by $\phi_{i,j}$ the $j^{th}$ column of the $i^{th}$ block $\Phi_i$, and write $Z_{0,ij} = \phi_{i,j}^*\Phi_\cl{S}(\Phi_\cl{S}^*\Phi_\cl{S})^{-1}\sign(\beta_\cl{S})$ for each $i \in \cl{S}^C, 1\le j \le m$. Additionally, denote $$Z_0 = \max_{i \notin \cl{S}, 1 \le j \le m}|Z_{0,ij}| = \|\Phi_{\cl{S}^C}^*\Phi_\cl{S}(\Phi_\cl{S}^*\Phi_\cl{S})^{-1}\sign(\beta_\cl{S})\|_\infty.$$ We simply need to show that with large probability $Z_0 \le 1/4$. For brevity, we write $W_{i,j} = (\Phi_\cl{S}^*\Phi_\cl{S})^{-1}\Phi_\cl{S}^*\phi_{i,j}$; we can then write $Z_{0,ij} = \langle W_{i,j},\sign(\beta)\rangle = \sum_{l \in \cl{S}}\langle W_{ij,l},\sign(\beta_l)\rangle$, where $W_{ij,l}$ denotes the block of the column $W_{i,j}$ corresponding to the entry $l \in \cl{S}$. We bound the magnitude of the sum terms as
\begin{align*}
|\langle W_{ij,l},\sign(\beta_l)\rangle| \le \|W_{ij,l}\|_2 \|\sign(\beta_l)\|_2 = \sqrt{m}\|W_{ij,l}\|_2.
\end{align*}
Next, we use Hoeffding's inequality and arguments similar to the ones preceding \eqref{eq:gamma1} in Appendix~\ref{sec:l12proof} to obtain that conditioned on a bound
\begin{align}
\gamma > \|\Phi_\cl{S}^*\Phi_{\cl{S}^C}\|_{B,1}
\label{eq:gamma}
\end{align}
and the invertibility condition, the probability of the inequality (\ref{eq:csc1}) failing to hold is at most $2pe^{-1/128m\gamma^2}$.

In order to establish \eqref{eq:csc2_1}, we use the following result that is a simple consequence of the Chernoff bound on the tail probability of the Gaussian distribution \cite{Kay.Book1998} and the union bound (see, e.g., \cite[Lemma 3.3]{CandesPlan}).
\begin{lem}
Let $(W_j')_{j \in J}$ be a fixed collection of vectors in $\real^n$ and set $Z_1 = \max_{j \in J} |\langle W_j',z\rangle|$. We then have $\prob(Z_1 \ge t) \le 2|J|e^{-t^2/2(\kappa')^2}$ for any $\kappa' \ge \max_{j \in J} \|W_j'\|_2$. \label{lem:randbound}
\end{lem}
\noindent We denote $W'_{i,j} = \Phi_\cl{S}(\Phi_\cl{S}^*\Phi_\cl{S})^{-1}\Phi_\cl{S}^*\phi_{i,j}$ for $i \notin \cl{S}$, $1 \le j \le m$, where $\phi_{i,j}$ represents the $j^{th}$ column of the $i^{th}$ block $\Phi_i$. Then we can write $$Z_1 = \|\Phi_{\cl{S}^C}^*\Phi_\cl{S}(\Phi_\cl{S}^*\Phi_\cl{S})^{-1}z\|_\infty = \max_{i \notin \cl{S}, 1\le j \le m}|\langle W'_{i,j},z\rangle|.$$ To use Lemma~\ref{lem:randbound} in this case, we again assume the invertibility condition holds and search for a bound on $\kappa'$:
\begin{align*}
\kappa' &= \max_{i \notin \cl{S},1 \le j \le m}\|\Phi_\cl{S}(\Phi_\cl{S}^*\Phi_\cl{S})^{-1}\Phi_\cl{S}^*\phi_{i,j}\|_2 \le \sqrt{6} \max_{i \notin \cl{S}, 1 \le j \le m} \|\Phi_\cl{S}^*\phi_{i,j}\|_2 \\
&\le \sqrt{6} \max_{i \notin \cl{S}} \|\Phi_\cl{S}^*\Phi_{i}\|_2 = \sqrt{6}\|\Phi_\cl{S}^*\Phi_{\cl{S}^C}\|_{B,1} \le \sqrt{6}\gamma.
\end{align*}
Thus, we have that conditioned on the bound (\ref{eq:gamma}) and the invertibility condition, (\ref{eq:csc2_1}) holds except with probability at most $2pe^{-(3/2-\sqrt{2})^2\lambda^2/4\gamma^2}$.

To finalize, we define $Z = \|\Phi_\cl{S}^*\Phi_\cl{S}-\ident\|_2$ and define the event $$E = \{Z \le 1/2\} \cup \{\|\Phi_\cl{S}^*\Phi_{\cl{S}^C}\|_{B,1} \le \gamma\}.$$ Then we have that the probability $P$ of the complementary size condition not being met is upper bounded by
\begin{align*}
P &\le \prob(\{Z_0 > 1/4\} \cup \{Z_1 \ge(3/2-\sqrt{2})\lambda\}|E)+\prob(E^C)\\
&\le 2pe^{-1/128m\gamma^2} + 2pe^{-(3/2-\sqrt{2})^2\lambda^2/4\gamma^2} + \prob(Z > 1/2) + \prob(\|\Phi_\cl{S}^*\Phi_{\cl{S}^C}\|_{B,1} > \gamma).
\end{align*}
We set $\gamma = c_3''/\sqrt{m\log p}$ with small enough $c_3''$ so that each of the first two terms of the right hand side is upper bounded by $2 p^{-4\log 2}$. To get the probability of the bound (\ref{eq:gamma}) not being valid, we once again appeal to the arguments made in relation to \eqref{eqn:Appendix.B.end} in Appendix~\ref{sec:l12proof} to get $\prob(\|\Phi_\cl{S}^*\Phi_{\cl{S}^C}\|_{B,1} > \frac{c_3''}{\sqrt{m\log{p}}}) < 2p^{-4 \log 2}$. Thus, the complementary size condition holds with probability at least $1-8p^{-4\log2}$.

By combining the three conditions (noting that the third condition already accounts for the first one), we have that Theorem~\ref{thm:lasso} holds with probability at least $1-8p^{-4\log2}-p^{-1}(2\pi \log p)^{-1/2}$.

\setcounter{lem}{0}
\section{Proof of Theorem~\ref{thm:grouplasso}}
\label{sec:grouplassoproof}
The proof of this theorem mirrors the steps taken for the proof of Theorem~\ref{thm:lasso} along with some necessary modifications. To begin, we need the following lemma concerning the behavior of the group lasso.
\begin{lem}
The group lasso estimate obeys $\|\Phi^*(y-\Phi\betahat)\|_{2,\infty} \le 2\lambda\sigma\sqrt{m}.$ \label{lem:gle}
\end{lem}
\begin{IEEEproof}
Since $\betahat$ minimizes the objective function over $\beta$, then 0 must be a subgradient of the objective function at $\betahat$. The subgradients of the group lasso objective function are of the form{\cite{YuanLin}} $$\Phi_i^*(\Phi\beta-y)+2\lambda\sigma\sqrt{m} \epsilon_i = 0,~i = 1,\ldots,r,$$ where $\epsilon_i \in \real^m$ is given by $\epsilon_i = \signbar(\beta_i)$ if $\beta_i \neq 0$ and $\|\epsilon_i\|_2 \le 1$ otherwise. Hence, since 0 is a subgradient at $\betahat$, there exists $\epsilon = [\epsilon_1^*~\ldots~\epsilon_r^*]^*$ such that $$\Phi^*(\Phi\betahat-y)=-2\lambda\sigma\sqrt{m} \epsilon.$$ The conclusion follows from the fact that $\|\epsilon\|_{2,\infty} \le 1$.
\end{IEEEproof}

In the following, we assume that $\sigma = 1$ without loss of generality and establish that the following three conditions together imply the theorem:
\begin{itemize}
\item {\em Invertibility condition:} The submatrix $\Phi_\cl{S}^*\Phi_\cl{S}$ is invertible and obeys $\|(\Phi_\cl{S}^*\Phi_\cl{S})^{-1}\|_2 \le 2.$
\item {\em Group orthogonality condition:} The vector $z$ satisfies the following inequality: $\|\Phi^*z\|_{2,\infty} \le \sqrt{2m}\cdot\lambda$.
\item {\em Group complementary size condition:} The following inequality holds: $$\|\Phi_{\cl{S}^C}^*\Phi_\cl{S}(\Phi_\cl{S}^*\Phi_\cl{S})^{-1}\Phi_\cl{S}^*z\|_{2,\infty} +2\lambda\sqrt{m}\|\Phi_{\cl{S}^C}^*\Phi_\cl{S}(\Phi_\cl{S}^*\Phi_\cl{S})^{-1}\signbar(\beta_\cl{S})\|_{2,\infty} \le (2-\sqrt{2})\lambda\sqrt{m}.$$
\end{itemize}
We assume that these three conditions hold. Since $\betahat$ minimizes the group lasso objective function, we must have $$\frac{1}{2}\|y-\Phi\betahat\|_2^2 + 2\lambda\sqrt{m}\|\betahat\|_{2,1} \le \frac{1}{2}\|y-\Phi\beta\|_2^2 + 2\lambda\sqrt{m}\|\beta\|_{2,1}.$$ Define $h := \betahat-\beta$, and note that $$\|y-\Phi\betahat\|_2^2 = \|(y-\Phi\beta)-\Phi h\|_2^2 = \|\Phi h\|_2^2+\|y-\Phi\beta\|_2^2-2\langle \Phi h,y-\Phi\beta\rangle.$$ Plugging this identity with $z=y-\Phi\beta$ into the above inequality and rearranging the terms gives $$\frac{1}{2}\|\Phi h\|_2^2 \le \langle \Phi h,z \rangle+2\lambda\sqrt{m}(\|\beta\|_{2,1}-\|\betahat\|_{2,1}).$$ Next, break up $\betahat$ into $\betahat_\cl{S}$ and $\betahat_{\cl{S}^C}=h_{\cl{S}^C}$ and rewrite the above equation as
\begin{align}
\frac{1}{2}\|\Phi h\|_2^2 \le \langle h,\Phi^*z \rangle+2\lambda\sqrt{m}(\|\beta_\cl{S}\|_{2,1}-\|\betahat_\cl{S}\|_{2,1} - \|h_{\cl{S}^C}\|_{2,1}). \label{eq:xh}
\end{align}
For each $i \in \cl{S}$, we have
\begin{align*}
\|\betahat_i\|_2 &= \|\beta_i+h_i\|_2 \ge \left|\left\langle \beta_i+h_i,\frac{\beta_i}{\|\beta_i\|_2}\right\rangle\right| \ge \left\langle \beta_i+h_i,\frac{\beta_i}{\|\beta_i\|_2}\right\rangle = \|\beta_i\|_2 + \langle h_i,\signbar(\beta_i)\rangle,
\end{align*}
where the first inequality is due to the projection of $\beta_i+h_i$ on $\mathrm{span}\{\beta_i\}$ having magnitude at most $\|\beta_i+h_i\|_2$. Thus, we can write $\|\betahat_\cl{S}\|_{2,1} \ge \|\beta_\cl{S}\|_{2,1}+\langle h_\cl{S},\signbar(\beta_\cl{S})\rangle$. Merging this inequality with~(\ref{eq:xh}) gives us
\begin{align}
\frac{1}{2}\|\Phi h\|_2^2 &\le \langle h,\Phi^*z \rangle-2\lambda\sqrt{m}(\langle h_\cl{S},\signbar(\beta_\cl{S})\rangle+\|h_{\cl{S}^C}\|_{2,1})\nonumber,\\
&= \langle h_\cl{S},\Phi_\cl{S}^*z\rangle + \langle h_{\cl{S}^C},\Phi_{\cl{S}^C}^*z\rangle-2\lambda\sqrt{m}(\langle h_\cl{S},\signbar(\beta_\cl{S})\rangle+\|h_{\cl{S}^C}\|_{2,1}).
\label{eq:xh2}
\end{align}
The group orthogonality condition and Lemma~\ref{lem:holder2} also imply $$\langle h_{\cl{S}^C},\Phi_{\cl{S}^C}^*z\rangle \le \|h_{\cl{S}^C}\|_{2,1}\|\Phi_{\cl{S}^C}^*z\|_{2,\infty} \le \sqrt{2m}\cdot\lambda\|h_{\cl{S}^C}\|_{2,1}.$$ Merging this result with (\ref{eq:xh2}) results in
\begin{align}
\frac{1}{2}\|\Phi h\|_2^2 &\le \langle h,\Phi^*z \rangle+2\lambda\sqrt{m}(-\langle h_\cl{S},\signbar(\beta_\cl{S})\rangle-\|h_{\cl{S}^C}\|_{2,1})\nonumber,\\
& \le |\langle h_\cl{S},v \rangle|-(2-\sqrt{2})\lambda\sqrt{m}\|h_{\cl{S}^C}\|_{2,1},
\label{eq:xh3}
\end{align}
where $v=\Phi_\cl{S}^*z-2\lambda\sqrt{m}\cdot\signbar(\beta_\cl{S})$. We aim to bound each of the terms on the right hand side independently. For the first term, we have
\begin{align*}
|\langle h_\cl{S},v\rangle| &= |\langle (\Phi_\cl{S}^*\Phi_\cl{S})^{-1}\Phi_\cl{S}^*\Phi_\cl{S}h_\cl{S},v\rangle| = |\langle \Phi_\cl{S}^*\Phi_\cl{S}h_\cl{S},(\Phi_\cl{S}^*\Phi_\cl{S})^{-1}v\rangle| \\
&\le |\langle \Phi_\cl{S}^*\Phi h,(\Phi_\cl{S}^*\Phi_\cl{S})^{-1}v\rangle|+|\langle \Phi_\cl{S}^*\Phi_{\cl{S}^C}h_{\cl{S}^C},(\Phi_\cl{S}^*\Phi_\cl{S})^{-1}v\rangle|.
\end{align*}
Denote the two terms on the right hand side as $A_1$ and $A_2$, respectively. For $A_1$ we use Lemma~\ref{lem:holder2} to obtain $$A_1 \le \|(\Phi_\cl{S}^*\Phi_\cl{S})^{-1}v\|_{2,1}\|\Phi_\cl{S}^*\Phi h\|_{2,\infty}.$$ Now we bound these two terms. For the first term, we get $$\|(\Phi_\cl{S}^*\Phi_\cl{S})^{-1}v\|_{2,1} \le \sqrt{k}\|(\Phi_\cl{S}^*\Phi_\cl{S})^{-1}v\|_2 \le \sqrt{k}\|(\Phi_\cl{S}^*\Phi_\cl{S})^{-1}\|_2\|v\|_2 \le 2k\|v\|_{2,\infty}$$ due to the invertibility condition. Using the group orthogonality condition, we also get $$\|v\|_{2,\infty} = \|\Phi_\cl{S}^*z-2\lambda\sqrt{m}\cdot\signbar(\beta_\cl{S})\|_{2,\infty} \le \|\Phi_\cl{S}^*z\|_{2,\infty}+2\lambda\sqrt{m}\le (2+\sqrt{2})\lambda\sqrt{m}.$$ For the second term $\|\Phi_\cl{S}^*\Phi h\|_{2,\infty}$, we use Lemma~\ref{lem:gle} and the group orthogonality condition to get
\begin{align*}
\|\Phi_\cl{S}^*\Phi h\|_{2,\infty} &\le \|\Phi_\cl{S}^*(\Phi\beta-y)\|_{2,\infty}+\|\Phi_\cl{S}^*(y-\Phi\betahat)\|_{2,\infty} \\
&= \|\Phi_\cl{S}^*z\|_{2,\infty}+\|\Phi_\cl{S}^*(y-\Phi\betahat)\|_{2,\infty} \le (2+\sqrt{2})\lambda\sqrt{m}.
\end{align*}
Combining, we finally get $A_1 \le 2(2+\sqrt{2})^2\lambda^2mk$. For $A_2$, we have from Lemma~\ref{lem:holder2} that $$A_2 \le \|h_{\cl{S}^C}\|_{2,1}\|\Phi_{\cl{S}^C}^*\Phi_\cl{S}(\Phi_\cl{S}^*\Phi_\cl{S})^{-1}v\|_{2,\infty} \le (2-\sqrt{2})\lambda\sqrt{m}\|h_{\cl{S}^C}\|_{2,1},$$ because of the group complementary size condition. Using now these bounds on $A_1,A_2$, we have $$|\langle h_\cl{S},v\rangle| \le 2(2+\sqrt{2})^2\lambda^2mk+(2-\sqrt{2})\lambda\sqrt{m}\|h_{\cl{S}^C}\|_{2,1}.$$ Plugging this into (\ref{eq:xh3}) gives $$\frac{1}{2}\|\Phi(\beta-\betahat)\|_2^2 \le 2(2+\sqrt{2})^2\lambda^2mk,$$ which suffices to prove the theorem, modulo the three conditions.

To finish the proof of the theorem, we now must evaluate the probability of each condition failing to hold under the assumed statistical model. The invertibility condition in this regard simply follows from Theorem~\ref{thm:rand_cond}, which means that it fails to hold with probability at most $2p^{-4\log 2}$. Next, note that $\|\Phi^*z\|_{2,\infty} \le \sqrt{2}\cdot\lambda\sqrt{m}$ is implied by $\|\Phi^*z\|_\infty \le \sqrt{2}\cdot\lambda$, which matches the orthogonality condition in the proof of Theorem~\ref{thm:lasso}  (cf. Appendix~\ref{app:thm:lasso}). Therefore, the group orthogonality condition fails to holds with probability at most $p^{-1}(2\pi \log p)^{-1/2}$ in the case of the group lasso. Therefore, we only need to evaluate the group complementary size condition.

In order to study the group complementary size condition, we partition it into two statements:
\begin{align}
Z_0 &:= \|\Phi_{\cl{S}^C}^*\Phi_\cl{S}(\Phi_\cl{S}^*\Phi_\cl{S})^{-1}\signbar(\beta_\cl{S})\|_{2,\infty} \le \frac{1}{4},\label{eq:csc3}\\
Z_1 &:= \|\Phi_{\cl{S}^C}^*\Phi_\cl{S}(\Phi_\cl{S}^*\Phi_\cl{S})^{-1}\Phi_\cl{S}^*z\|_{2,\infty} \le \left(\frac{3}{2}-\sqrt{2}\right)\lambda\sqrt{m}.\label{eq:csc4}
\end{align}
In order to evaluate \eqref{eq:csc3}, we compare it to (\ref{eq:l12cond}) and note that the main difference between the two expressions is a change from $1/4$ to $1$. Given that both Theorem~\ref{thm:l12} and this theorem operate under the same statistical model, it is therefore straightforward to argue from the analysis of (\ref{eq:l12cond}) that (\ref{eq:csc3}) holds except with probability at most $2pe^{-1/128m\gamma^2}$, where $\gamma$ is defined as any positive scalar that satisfies $\gamma \geq \|\Phi_\cl{S}^*\Phi_{\cl{S}^C}\|_{B,1}$. The second condition (\ref{eq:csc4}) is implied by the inequality
\begin{align}
Z_2 :=  \|\Phi_{\cl{S}^C}^*\Phi_\cl{S}(\Phi_\cl{S}^*\Phi_\cl{S})^{-1}\Phi_\cl{S}^*z\|_\infty &\le\left(\frac{3}{2}-\sqrt{2}\right)\lambda,
\label{eq:csc2}
\end{align}
which is shown to hold with probability at most $2pe^{-(3/2-\sqrt{2})^2\lambda^2/4\gamma^2}$ in (\ref{eq:csc2_1}) (cf.~Appendix~\ref{app:thm:lasso}). To conclude, we once again define $Z = \|\Phi_\cl{S}^*\Phi_\cl{S}-\ident\|_2$ and define the event $$E = \{Z \le 1/2\} \cup \{\|\Phi_\cl{S}^*\Phi_{\cl{S}^C}\|_{B,1} \le \gamma\}.$$ Then we have that the probability $P$ of the group complementary size condition not being met is upper bounded by
\begin{align*}
P &\le \prob(\{Z_0 > 1/4\} \cup \{Z_1 \ge(3/2-\sqrt{2})\lambda\}|E)+\prob(E^C)\\
&\le 2pe^{-1/128m\gamma^2} + 2pe^{-(3/2-\sqrt{2})^2\lambda^2/12\gamma^2} + \prob(Z > 1/2) + \prob(\|\Phi_\cl{S}^*\Phi_{\cl{S}^C}\|_{B,1} > \gamma)\\
&\le 2pe^{-1/128m\gamma^2} + 2pe^{-(3/2-\sqrt{2})^2\lambda^2/12\gamma^2} + 2p^{-4\log 2} +\prob(\|\Phi_\cl{S}^*\Phi_{\cl{S}^C}\|_{B,1} > \gamma).
\end{align*}
We notice that this inequality has been shown to be upper bounded by $8p^{-4\log2}$ in Appendix~\ref{app:thm:lasso}. Thus, the group complementary size condition fails to hold with probability at most $8p^{-4\log2}$.

By combining the failures of the three conditions (and noting that the third condition already accounts for the first one), we have that Theorem~\ref{thm:grouplasso} holds with probability at least $1-8p^{-4\log2}-p^{-1}(2\pi \log p)^{-1/2}$.

\ifCLASSOPTIONcaptionsoff
  \newpage
\fi


\end{document}